\documentclass[11pt,draftcls,onecolumn]{IEEEtran}
\usepackage{mathptmx}
\usepackage [T1]{fontenc}
\usepackage{amsmath}
\usepackage{amssymb}
\usepackage{amsfonts}
\usepackage{textcomp}
\usepackage{graphicx}
\usepackage{rotating}
\usepackage{mathdots}
\usepackage{color}
\usepackage{bm}

\definecolor{purple}{rgb}{0.7,0,0.4}
\definecolor{darkgreen}{rgb}{0,0.5,0}

\newcommand{\matc}[1]{\ensuremath{{\mathcal{#1}}}}
\newcommand{\mat}[1]{\ensuremath{{\mathbf{#1}}}}

\newcommand{\inmat}[3][C]{\ensuremath{\in{\mathbb{#1}}^{{#2}\times{#3}}}}
\newcommand{\inten}[4][C]{\ensuremath{\in{\mathbb{#1}}^{{#2}\times{#3}\times{#4}}}}
\newcommand{\be}{\begin{equation}}
\newcommand{\ee}{\end{equation}}
\newcommand{\beq}{\begin{eqnarray}}
\newcommand{\eeq}{\end{eqnarray}}
\newcommand{\bi}{\begin{itemize}}
\newcommand{\ei}{\end{itemize}}
\newcommand{\ba}{\begin{array}}
\newcommand{\ea}{\end{array}}

\begin{document}
\title{Consensus-based In-Network Computation of the PARAFAC Decomposition}

\author{Alain Y. Kibangou$^{*}$ and André L. F. de Almeida
\thanks{\noindent
*Alain Y. Kibangou is with GIPSA-Lab, University Joseph Fourier, CNRS, 11 rue des Mathématiques, Grenoble Campus, 38402 Saint Martin d'Hères Cedex, France (Phone: +33 4 76 82 64 51, e-mail: alain.kibangou@ujf-grenoble.fr). Andr{\'e} L. F. de Almeida is with the Department of Teleinformatics Engineering, Federal University of Cear{\'a}, CP 6005, 60455-760, Fortaleza, Brazil (Phone/Fax: +55 85 3366 9470, e-mail: andre@gtel.ufc.br. This work is partially supported by FUNCAP and INRIA. Andr{\'e} L. F. de Almeida is partially supported by the CNPq.}}

\maketitle

\vspace{-5ex}
\begin{abstract}
Higher-order tensor analysis is a multi-disciplinary tool widely used in numerous application areas involving data analysis such as psychometrics, chemometrics, and signal processing, just to mention a few. The parallel factor (PARAFAC) decomposition, also known by the acronym CP (standing for ``CANDECOMP/PARAFAC'' or yet ``canonical polyadic'') is the most popular tensor decomposition. Its widespread use comes from its essential uniqueness property under mild conditions as well as to the existence of several numerical algorithms that can be used to compute the decomposition. 
In this work, we present a new approach for the distributed computation of the PARAFAC decomposition of a third-order tensor across a network of collaborating nodes. We are interested in the case where the overall data gathered across the network can be modeled as a data tensor admitting an essentially unique PARAFAC decomposition, while each node only observes a sub-tensor with not necessarily enough diversity so that identifiability conditions are not locally fulfilled at each node. In this situation, conventional (centralized) tensor based methods cannot be applied individually at each node.
By allowing collaboration between neighboring nodes of the network, we propose distributed versions of the alternating least squares (ALS) and Levenberg-Marquardt (LM) algorithms for the in-network estimation of the factor matrices of a third-order tensor. We assume that one of the factor matrices contains parameters that are local to each node, while the two remaining factor matrices contain global parameters that are common to the whole network. The proposed algorithms combine the estimation of the local factors with an in-network computation of the global factors of the PARAFAC decomposition using average consensus over graphs. They emulate their centralized counterparts in the case of ideal data exchange and ideal consensus computations. The performance of the proposed algorithms are evaluated in both ideal and imperfect cases.
\end{abstract}
\vspace{-2ex}
\begin{keywords}
\vspace{-1ex}
\noindent Distributed estimation, collaborative networks, PARAFAC decomposition, average consensus.
\end{keywords}



\section{Introduction}
\label{sec:intro}
In a various number of application domains, observed data exhibit intrinsically a multiway structure. The samples can be indexed by several indices, which can be directly related to a physical meaning. The more obvious data or signals with such a multiway structure are images. However, during the last decade, several applications in signal processing for communications were devoted to other kind of multidimensional signals \cite{SiGiBr:00,AlFaMo:07,NiDL:08,KiFa:09e}. For higher-order multiway arrays, also called tensors, i.e. extensions of the notion of matrix (two-way array) to greater order, the use of multilinear algebra tools is highly recommended in order to fully exploit the inherent structure of these signals. These tools have been used for both analysis and synthesis of existing and new communications systems for instance. The most used tensor tools are certainly the Tucker decomposition \cite{Tuc:66}, or higher-order singular value decomposition \cite{DeDeVa:00}, and the parallel factor (PARAFAC)  decomposition \cite{Har:70}, also known by the acronym CP, which stands for ``CANDECOMP/PARAFAC'' \cite{CaCh:70}  or, alternatively, ``canonical polyadic'' decomposition \cite{ComoLA5:JC}.

PARAFAC consists in decomposing an $N$-order tensor in a sum of rank-one tensors. Its widespread use is due to its essential uniqueness property under mild conditions \cite{Kru:77,SiBr:00,StSi:07} as well as to the existence of several numerical algorithms that can be used to compute this decomposition \cite{ToBr:06,ComoLA5:JC}. The uniqueness property is a key issue when trying to infer some information from the PARAFAC decomposition of a given tensor.

In most of applications in the literature, batch processing is considered for computing the PARAFAC decomposition. Several PARAFAC fitting algorithms have been proposed in the literature. They can be classified into three main categories: alternating algorithms, derivative based algorithms, and non-iterative algorithms (see \cite{ToBr:06} for a comparison of various algorithms). In alternating algorithms, only a subset of the parameters are updated at each step whereas in derivative based algorithms all the parameters are simultaneously updated by successive approximations. To the first category belongs the alternating least squares (ALS) algorithm \cite{SmBrGe:04} while the damped Gauss-Newton algorithm, also known as Levenberg-Marquardt (LM) algorithm, belong to the second category \cite{ToBr:06}. Recently, by taking into account the fact that data can be serially acquired or the underlying process can be time-varying, adaptive algorithms have been proposed in \cite{NiSi:09}. These algorithms are initialized by using a batch processing.

In this paper, motivated by communication issues in wireless sensor networks, we consider the problem of in-network computation of the PARAFAC decomposition across a network of communicating nodes. Indeed, in recent years, the deployment of sensors for monitoring, collaborative information processing and control has drawn a considerable attention. In particular, wireless sensor networks that can operate autonomously, i.e. without a fusion center collecting and processing all measurements, exhibit desirable properties such as robustness against node failure \cite{OSFaMu:07}. The coordinated action of different nodes requires local exchange of information.

In order to exploit multiple forms of diversity present in communication signals, one can resort to the powerful tensor-based blind estimation methods proposed in the literature for multiuser communications (see, e.g. \cite{SiGiBr:00,SiDi:01,AlFaMo:07,NiDL:08,KiFa:09e}) and, more recently, for multiple-input multiple-output (MIMO) communications (see \cite{SiBu:02,AlFaMo:08c,AlFaMo:08b,AlFaMo:08d} and references therein). However, these methods assume that a given node can compute the PARAFAC decomposition of the observed data tensor, which is generally possible by exploiting multiple dimensions of the observed data such as space, time, frequency and/or code. However, many wireless devices are limited by size, hardware complexity, or other constraints such as the use of a single (or a very small number of) antenna element(s). Hence, even though the overall observations gathered across the network of nodes can be modeled as a data tensor with enough diversity for information recovery purposes, each node only has access to a sub-tensor with not necessarily enough diversity. In this case, standard tensor based methods cannot be applied for such nodes. Indeed, without collaboration with its neighbors, the nodes cannot retrieve the factors of the PARAFAC decomposition.

In an attempt to filling this gap, the authors have derived a distributed alternating least squares algorithm for joint channel and symbol estimation in DS-CDMA networks \cite{KiAl:10}. This algorithm has also been applied to distributed decoding in distributed space-time block coded systems \cite{KiAl:11}. The approach proposed in these works are based on the concept of average consensus in graphs, which is an important issue in algorithm design for distributed computing. Average consensus has been extensively studied in computer science (distributed agreement and synchronization problems for example) and is a central topic for load balancing (with divisible tasks) in parallel computers. It has recently found application in distributed coordination of mobile autonomous agents, distributed data fusion in sensor networks, and distributed estimation and control \cite{XiBoKi:07}.

Following the idea of \cite{KiAl:10}, this work presents a generalized approach for distributed computation of the PARAFAC decomposition of a third-order tensor in a collaborative network. We assume that each node has access to several slices of the overall third-order tensor, hence a third-order sub-tensor. We are interested in the case where the overall data gathered across the network can be modeled as a data tensor admitting an essentially unique PARAFAC decomposition, while each node only observes a sub-tensor with not necessarily enough diversity so that identifiability conditions are not locally fulfilled. In such a situation, conventional (centralized) tensor based methods cannot be applied individually at each node, unless some form of collaboration is allowed in the network. 
By resorting to the concept of average consensus in graphs, we propose distributed versions of the alternating least squares (ALS) and Levenberg-Marquardt (LM) algorithms for the in-network estimation of the factor matrices of a PARAFAC decomposition of a third-order tensor. We assume that one of the factor matrices contains parameters that are local to each node, while the two remaining factor matrices contain global parameters that are common to the whole network. The proposed algorithms combine the estimation of the local factors with an in-network computation of the global factors of the PARAFAC decomposition using average consensus iterations.

In summary, the contribution of this work is two-fold. First, the algorithm development presented here provides a generalized framework by considering that collaborating nodes observe third-order tensors, whose dimensionality may differ in a particular mode. Second, in contrast to \cite{KiAl:10} which was restricted to ALS fitting, we develop a distributed version of LM algorithm that presents a faster convergence. With perfect average consensus, our numerical results show that these distributed PARAFAC-based algorithms yield similar performances as their centralized counterparts. We also evaluate the impact of noisy exchanges and imperfections in average consensus computations. Without loss of generality, we restrict our study to third-order tensors.

The paper is organized as follows. In Section \ref{sec_problem_statement}, the problem is formulated along with preliminaries including the PARAFAC decomposition and consensus problems in graphs. In Sections \ref{sec_DALS} and \ref{sec_DLM}, we derive distributed algorithms for computing the PARAFAC decomposition. We first propose a distributed alternating least squares (ALS) algorithm in Section \ref{sec_DALS}, and then a distributed Levenberg-Marquardt (LM) one in Section \ref{sec_DLM}. The efficiency of the proposed algorithms are evaluated in Section \ref{sec:simres} by means of computer simulations, before concluding the paper in Section \ref{sec:conc} .

\par\noindent\textit{Notations}: Vectors are written as boldface lower-case letters ($\mathbf{a}$,$\mathbf{b}$,$\cdots$), matrices as boldface upper-case letters ($\mathbf{A}$,$\mathbf{B}$,$\cdots$), and tensor as blackboard letters ($\mathbb{A}$,$\mathbb{B}$,$\cdots$). $\mathbf{A}_{i.}$ and $\mathbf{A}_{.j}$ denote respectively the $i$th row and the $j$th column of the $I\times J$ matrix $\mathbf{A}$. $\mathbf{A}^T$ stands for the transpose of $\mathbf{A}$ whereas $\mathbf{A}^H$ stands for its complex conjugate. $diag(.)$ is the operator that forms a diagonal matrix from its vector argument whereas $vec(.)$ forms a vector by stacking the columns of its matrix argument. For matrices $\mathbf{X}$ and $\mathbf{Y}$ with the same dimension, $\star$ stands for the Schur-Hadamard matrix product, i.e. the entrywise product. For $\mathbf{X}\in\mathbb{C}^{I\times R}$, and $\mathbf{Y}\in\mathbb{C}^{J\times R}$, the Khatri-Rao product,
denoted by $\odot$, is defined as follows:
\begin{equation}\label{khatri_rao}
\mathbf{X}\odot\mathbf{Y}=\left(
                             \begin{array}{ccc}
                               \mathbf{Y}diag(\mathbf{X}_{1.})\\
                               \mathbf{Y}diag(\mathbf{X}_{2.})\\
                               \vdots\\
                               \mathbf{Y}diag(\mathbf{X}_{I.})\\
                             \end{array}
                           \right)\in\mathbb{C}^{IJ\times R}.
\end{equation} In the sequel, we will make use of the properties of the Khatri-Rao product
and the Frobenius norm given below:
\begin{eqnarray}\label{prop_khatri2} vec(\mathbf{X} diag(\mathbf{z})
\mathbf{Y}^T) &= &(\mathbf{Y}\odot\mathbf{X})\mathbf{z},\\
\label{prop_khatri}
\mathbf{X}\odot\mathbf{Y}&=&\pmb{\Pi}(\mathbf{Y}\odot\mathbf{X}),\\
\label{prop_khatri_block}
\left(
  \begin{array}{c}
    \mathbf{X}_1 \\
    \vdots \\
    \mathbf{X}_M \\
  \end{array}
\right)\odot\mathbf{Y}&=&\left(
  \begin{array}{c}
    \mathbf{X}_1\odot\mathbf{Y} \\
    \vdots \\
    \mathbf{X}_M\odot\mathbf{Y} \\
  \end{array}
\right),\\
\label{prop_frobenius}
\left\|\mathbf{X}\right\|_F^2&=&\left\|\pmb{\Pi}\mathbf{X}\right\|_F^2,\\
\label{prop_frobenius2}
\left\|\left(
         \begin{array}{ccc}
           \mathbf{X}_1 & \cdots & \mathbf{X}_M \\
         \end{array}
       \right)
\right\|_F^2&=&\sum\limits_{m=1}^M\left\|\mathbf{X}_m\right\|_F^2,\\
\label{prop_frobenius3}
\left\|\left(
         \begin{array}{ccc}
           \mathbf{x}_1 & \cdots & \mathbf{x}_M \\
         \end{array}
       \right)
\right\|_F^2&=&\sum\limits_{m=1}^M\left\|\mathbf{x}_m\right\|_2^2,
\end{eqnarray}
where $\pmb{\Pi}$, a permutation matrix, and $\mathbf{X}$ are matrices
with compatible dimensions. We also have:
\begin{equation}\label{product_KR}
\left(\mathbf{A}\odot\mathbf{B}\right)\left(\mathbf{C}\odot\mathbf{D}\right)=\mathbf{AC}\star\mathbf{BD},
\end{equation}
for matrices with compatible dimensions.
\section{Preliminaries and Problem Statement}
\label{sec_problem_statement} Let us consider a network of $L$ nodes whose communication links are modeled by means of an undirected graph $\mathcal{G}\left(\mathcal{N},\mathcal{E}\right)$. $\mathcal{N}=\left\{1,\cdots,L\right\}$ and $\mathcal{E}$ denote respectively the node set and the edge set, where each edge $(l,l')\in\mathcal{E}$ is an unordered pair of distinct nodes. We assume that $\mathcal{G}$ is a connected graph, meaning that it exists a path between any two pair of nodes. In other words, information can be exchanged by any pair of nodes through direct links or multi-hop ones.

Let us assume that the node $l\in\mathcal{N}$ has at its disposal the data $x_{i,j,k}^{(l)}$, $i=1,\ldots,I_l$, $j=1,\cdots,J$, $k=1,\cdots,K$. Assume that these data can be organized in a $I_l\times J\times K$ tensor $\mathbb{X}^{(l)}$ admitting a PARAFAC decomposition \cite{Har:70,CaCh:70}:
\begin{equation}\label{parafac_model}
x_{i,j,k}^{(l)}=\sum\limits_{r=1}^{R}a_{i,r}^{(l)}b_{j,r}c_{k,r}.
\end{equation}
The tensor $\mathbb{X}^{(l)}=[x_{i,j,k}^{(l)}]$ is completely characterized by three factor matrices $\mathbf{A}^{(l)}=[a_{i,r}^{(l)}]\in\mathbb{C}^{I_l\times R}$, $\mathbf{B}=[b_{j,r}]\in\mathbb{C}^{J\times R}$, and $\mathbf{C}=[c_{k,r}]\in\mathbb{C}^{K\times R}$. The purpose of the PARAFAC decomposition at each node is to find these factor matrices from the observed data tensor $\mathbb{X}^{(l)}$, $l=1,\ldots, L$.

\subsection{The PARAFAC decomposition}
PARAFAC is certainly the most popular tensor model that can be found in the literature. That is surely due to its essential uniqueness, which means that each factor matrix can be determined up to column scaling and permutation, i.e. two sets of matrices $\{\mathbf{A}^{(l)},\mathbf{B},\mathbf{C}\}$ and $\{\tilde{\mathbf{A}}^{(l)},\tilde{\mathbf{B}},\tilde{\mathbf{C}}\}$ giving rise to the same tensor $\mathbb{X}^{(l)}$ are linked by the following relations $\tilde{\mathbf{A}}^{(l)}=\mathbf{A}^{(l)}\pmb{\Pi\Delta}_A$, $\tilde{\mathbf{B}}=\mathbf{B}\pmb{\Pi\Delta}_B$, $\tilde{\mathbf{C}}=\mathbf{C}\pmb{\Pi\Delta}_C$, with $\pmb{\Delta}_A\pmb{\Delta}_B\pmb{\Delta}_C=\mathbf{I}_{R}$, where $\pmb{\Pi}$ is a $R\times R$ permutation matrix, whereas $\pmb{\Delta}_A$,$\pmb{\Delta}_B$, and $\pmb{\Delta}_C$ are $R\times R$ diagonal matrices. A sufficient condition for such an uniqueness, the so-called  Kruskal's condition, states that the PARAFAC decomposition  (\ref{parafac_model}) is essentially unique if \cite{Kru:77,SiGiBr:00}
\begin{equation}\label{Kruskal_generalized}
k_{\mathbf{A}^{(l)}}+k_{\mathbf{B}}+k_{\mathbf{C}}\geq 2R+2,
\end{equation}
where $k_{\mathbf{X}}$ denotes the Kruskal-rank, also called k-rank, of a given matrix $\mathbf{X}$, and corresponds to the greatest integer $k_{\mathbf{X}}$ such that any set of $k_\mathbf{X}$ columns of $\mathbf{X}$ is independent. The rank and the Kruskal-rank of $\mathbf{X}$ are linked by the following inequality $k_{\mathbf{X}}\leq rank(\mathbf{X})$.

Another feature of PARAFAC is to provide a simple link between the unfolded forms of a tensor and its factor matrices. The unfolded forms are obtained by concatenating tensor slices along the same mode. We call ``slice'', the matrix obtained by fixing one (out of the three) dimension of the tensor. For the third-order tensor $\mathbb{X}^{(l)}$, we have the following three types of slices:

\[\mathbf{X}_{i..}^{(l)}=\left(\begin{array}{ccc}
                         x_{i,1,1}^{(l)} & \cdots & x_{i,1,N}^{(l)} \\
                         \vdots & \ddots & \vdots \\
                         x_{i,J,1}^{(l)} & \cdots & x_{i,J,N}^{(l)} \\
                         \end{array}
                         \right)=\mathbf{B}diag(\mathbf{A}_{i.}^{(l)})\mathbf{C}^T,
\]
\[\mathbf{X}_{.j.}^{(l)}=\left(\begin{array}{ccc}
                         x_{1,j,1}^{(l)} & \cdots & x_{I_l,j,1}^{(l)} \\
                         \vdots & \ddots & \vdots \\
                         x_{1,j,K}^{(l)} & \cdots & x_{I_l,j,K}^{(l)} \\
                         \end{array}
                         \right)=\mathbf{C}diag(\mathbf{B}_{j.})\mathbf{A}^{(l)T},
\]
and
\[
\mathbf{X}_{..k}^{(l)}=\left(\begin{array}{ccc}
                         x_{1,1,k}^{(l)} & \cdots & x_{I_l,J,k}^{(l)} \\
                         \vdots & \ddots & \vdots \\
                         x_{1,1,k}^{(l)} & \cdots & x_{I_l,J,k}^{(l)} \\
                         \end{array}
                         \right)=\mathbf{A}^{(l)}diag(\mathbf{C}_{k.})\mathbf{B}^T.
\]
Using (\ref{khatri_rao}), the three vertically unfolded matrix
representations for $\mathbb{X}^{(l)}$ are respectively given by:
\be\label{eq:unfold_X1l}
\mathbf{X}^{(l)}_{1}=\left(
               \begin{array}{c}
                 \mathbf{X}_{.1.}^{(l)} \\
                 \vdots \\
                 \mathbf{X}_{.J.}^{(l)} \\
               \end{array}
             \right)=\left(\mathbf{B}\odot\mathbf{C}\right)\mathbf{A}^{(l)T}\in\mathbb{C}^{JK\times I_l},
\ee
\be\label{eq:unfold_X2l}
\mathbf{X}^{(l)}_{2}=\left(
               \begin{array}{c}
                 \mathbf{X}_{..1}^{(l)} \\
                 \vdots \\
                 \mathbf{X}_{..K}^{(l)} \\
               \end{array}
             \right)=\left(\mathbf{C}\odot\mathbf{A}^{(l)}\right)\mathbf{B}^T\in\mathbb{C}^{KI_l\times J},
\ee
and
\be\label{eq:unfold_X3l}
\mathbf{X}^{(l)}_{3}=\left(
               \begin{array}{c}
                 \mathbf{X}_{1..}^{(l)} \\
                 \vdots \\
                 \mathbf{X}_{I_l..}^{(l)} \\
               \end{array}
             \right)=\left(\mathbf{A}^{(l)}\odot\mathbf{B}\right)\mathbf{C}^T\in\mathbb{C}^{I_lJ\times K}.
\ee
Note that if $\mathbf{A}^{(l)}$ is a row vector, meaning that the measurements available at node $l$ are two-dimensional, we have the following equalities:
\begin{eqnarray}\label{particular_case}
&&\mathbf{X}^{(l)}_{3}=\mathbf{B}diag(\mathbf{A}^{(l)})\mathbf{C}^T=(\mathbf{X}^{(l)}_{2})^T\\
&&\textit{vec}(\mathbf{X}^{(l)T}_{3})=\mathbf{X}^{(l)}_{1}.
\end{eqnarray}

\subsection{Problem statement}
Each node can use  any fitting algorithm for computing its own PARAFAC decomposition. As stated above, essential uniqueness is guaranteed if the Kruskal's condition (\ref{Kruskal_generalized}) is fulfilled at each node. Herein, we assume that this condition is violated at each node, i.e $k_{\mathbf{A}^{(l)}}+k_{\mathbf{B}}+k_{\mathbf{C}}<2R+2, \quad\forall l\in\mathcal{N}$, and there is no local uniqueness property.

 If all measurements $x_{i,j,k}^{(l)}$, $l=1,\cdots,L$, $i=1,\ldots,I_l$, $j=1,\cdots,J$, $k=1,\cdots,K$, were available at a given central point (e.g. a fusion center), then one could define a global tensor $\mathbb{X} \inten{I}{J}{K}$, with $I=\sum\limits_{l=1}^LI_l$, which concatenates the sub-tensors $\mathbb{X}^{(l)}$, $l=1,\ldots, L$, along its first mode$\footnote{Without loss of generality, we have concatenated the $L$ sub-tensors along the first mode, since the PARAFAC decomposition of each sub-tensor differ in the first mode factor matrix. However, concatenation along the second or third modes would be possible if the $L$ PARAFAC decompositions differed in their second or third modes, respectively. In any case, the approach developed in this paper supposes that the PARAFAC decomposition of the $L$ subtensors differ only in one mode, the other ones being comon to the $L$ nodes.}$, i.e. $\mathbb{X}~=[\mathbb{X}^{(1)} \sqcup_1 \mathbb{X}^{(1)} \sqcup_1 \cdots \sqcup_1 \mathbb{X}^{(L)}]$. Let $x_{i,j,k}$ be the typical element of the global tensor $\mathbb{X}$. Note that
$$x_{i,j,k}=x^{(l)}_{i',j,k}, \quad \textrm{with} \quad i=\sum\limits_{m=1}^{l-1}I_m+i',$$
$i'=1,\ldots, I_l$, $l=1,\ldots, L$. Obviously, such a global tensor admits a PARAFAC model with $\mathbf{A}$, $\mathbf{B}$, and $\mathbf{C}$ as factor matrices, where $\mathbf{A}\in\mathbb{C}^{I\times R}$ results in a row-wise concatenation of $\mathbf{A}^{(l)}$, $l=1,\cdots,L$, i.e. $\mat A=[\mat A^{(1)T}, \ldots, \mat A^{(L)T}]^T$. The global estimation problem consists in minimizing the following cost function:
\beq\label{eq:global1}
\mathcal{J}(\mat A,\mat S,\mat C)&=& \frac{1}{2}\sum\limits_{i=1}^{I}\sum\limits_{j=1}^{J}
\sum\limits_{k=1}^{K}\left| x_{i,j,k} - \sum\limits_{r=1}^{R}a_{i,r}b_{j,r}c_{k,r}\right|^2 \nonumber\\
&=&\frac{1}{2}\sum\limits_{l=1}^{L}\sum\limits_{i=1}^{I_l}\sum\limits_{j=1}^{J}
\sum\limits_{k=1}^{K}\left| x_{i,j,k}^{(l)} - \sum\limits_{r=1}^{R}a_{i,r}^{(l)}b_{j,r}c_{k,r}\right|^2.
\eeq
We assume that the global uniqueness condition is fulfilled, i.e.
\[
k_{\mathbf{A}}+k_{\mathbf{C}}+k_{\mathbf{S}}\geq 2R+2.
\]

\noindent Note that unfolded matrices of the (global) tensor $\mathbb{X}$ can be built from those of the (local) sub-tensors $\mathbb{X}^{(l)}$ as follows:
\begin{equation}\label{unfolded_1}
\mathbf{X}_1=\left(
               \begin{array}{ccc}
                 \mathbf{X}^{(1)}_{1}&\cdots & \mathbf{X}^{(L)}_{1}
               \end{array}
             \right)=(\mathbf{B}\odot\mathbf{C})\mathbf{A}^T,
\end{equation}

\begin{equation}\label{unfolded_2}
\mathbf{X}_2=\pmb{\Pi}\left(
               \begin{array}{c}
                 \mathbf{X}^{(1)}_{2} \\
                 \vdots \\
                 \mathbf{X}^{(L)}_{2} \\
               \end{array}
             \right)=(\mathbf{C}\odot\mathbf{A})\mathbf{B}^T,
\end{equation}
\begin{equation}\label{unfolded_3}
\mathbf{X}_3=\left(
               \begin{array}{c}
                 \mathbf{X}^{(1)}_{3} \\
                 \vdots \\
                 \mathbf{X}^{(L)}_{3} \\
               \end{array}
             \right)=(\mathbf{A}\odot\mathbf{B})\mathbf{C}^T,
\end{equation}
where $\pmb{\Pi}$ stands for an appropriate permutation matrix.

Uniqueness of PARAFAC can be exploited by sending the sub-tensors $\mathbb{X}^{(1)}, \ldots, \mathbb{X}^{(L)}$ to a central node, or fusion center, from which the global tensor $\mathbb{X}$ can be constructed. The central node can then perform the PARAFAC decomposition of $\mathbb{X}$ and send the estimated factor matrices back to the nodes. However, it is well known that the existence of a central node in a network gives rise to a particularly vulnerable setup. Resorting to distributed estimation is then well suited. One could imagine that nodes exchange their received data samples with their neighbors. As a consequence, after such an information exchange, from its own data matrix and those received from its neighbors, each node can built a tensor, which is in fact a sub-tensor of $\mathbb{X}$. Unfortunately, we cannot ensure that all sub-tensors inherit the uniqueness property of the global tensor. Therefore, in such a scheme, some nodes can obtain undesirable estimates, i.e. estimates that cannot be linked to the actual factor matrices in a unique way. The purpose of the following section is to derive distributed estimation schemes, preserving the PARAFAC uniqueness property. For this purpose, we will resort to the notion of average consensus in graphs.

\subsection{Average consensus algorithm}
Average consensus is an important problem in algorithm design for distributed computing. It aims at reaching an agreement on a common value across a network; such a common value being the average of the initial local values.

Let $\mathbf{R}_l(0)$ be a matrix assigned to node $l$ at time $t=0$. The distributed average consensus problem consists in computing the average $(1/L)\sum\limits_{l=1}^{L}\mathbf{R}_l(0)$ at every node, via local communication and computation across the network modeled with an undirected graph $\matc{G}(\matc{N},\matc{E})$. Node $l$ carries out its update, at each step, based on its local state and communication with its neighboring nodes $l'\in\matc N_l=\big\{l'|(l,l') \in \matc E\big\}$.

There are several simple methods for distributed average consensus. For example, each node can store a table of all initial node values known at that time. At each step each pair of neighbors exchange tables of initial values and update their tables. In this flooding algorithm, all nodes know all initial values in a number of steps equal to the diameter of the graph, at which point each of them can compute the average \cite{XiBoKi:07}. In widely used average consensus algorithms, each node updates itself by adding a weighted sum of local discrepancies, i.e. differences between neighboring node values and its own value. In matrix form, we get:
\begin{equation}\label{consensus}
\mathbf{R}_l(t+1)=\mathbf{R}_l(t)+\sum\limits_{l'\in\mathcal{N}_l}w_{l,l'}\left(\mathbf{R}_{l'}(t)-\mathbf{R}_l(t)\right),
\end{equation}
where $w_{l,l'}$ is a weight associated with the edge $\{l,l'\}$, which is a typical element of the weight matrix $\mat W \inmat{L}{L}$, which is doubly stochastic, i.e. $\mat 1^T\mat W=\mat 1^T$ and $\mat W\mat 1= \mat 1$, where $\mat 1 \inmat{L}{1}$ is a vector of ones. If $1$ is a single eigenvalue of $\mathbf{W}$ and the remaining eigenvalues have magnitude lower than 1 then asymptotic convergence is guaranteed, i.e. $\mathbf{R}_l(t)\rightarrow (1/L)\sum\limits_{l=1}^{L}\mathbf{R}_l(0)$ as $t$ goes towards infinity. Several weights fulfilling the previous conditions have been proposed in the literature. One of them is given by $\mathbf{W}$ as $\mathbf{W}=\mathbf{I}-\gamma\mathbf{L}$ where $\mathbf{L}$ denotes the Laplacian of the graph whose entries $l_{pq}$ are given by:
\[
l_{pq}=\left\{\begin{array}{cc}
                d_p & \mathrm{if}\:\: p=q \\
                -1 & \mathrm{if}\:\: q\in\mathcal{N}_p \\
                0 & \mathrm{elsewhere}.
              \end{array}
\right.
\] In order to ensure convergence $\gamma$ should be selected as $\gamma <\frac{2}{\max\{\lambda_p\}}$, $p=1,\ldots, L$, $\{\lambda_p\}$ being the set of eigenvalues of the Laplacian matrix ordered as $0=\lambda_1<\lambda_2<\cdots<\lambda_L$. Since $\gamma$ modifies the spectrum of $\mat W$ and in particular its second largest eigenvalue, the choice of this parameter influences the speed of convergence of the consensus algorithm. It has been shown that the fastest policy is obtained by selecting $\gamma=\frac{2}{\lambda_2+\lambda_L}$.

Recently, \cite{Kib:11cdc,Kib:12acc} have shown that, in the noiseless case, average consensus can be achieved in a finite number of steps $D$ corresponding to the number of nonzero distinct eigenvalues of the Laplacian matrix. Therefore, provided the information on the spectrum of the Laplacian is available, the finite-time average consensus algorithm can be run as follows:
\begin{equation}\label{ft_consensus}
\mathbf{R}_l(t+1)=\mathbf{R}_l(t)+\frac{1}{\lambda_{p(t)}}\sum\limits_{l'\in\mathcal{N}_l}\left(\mathbf{R}_{l'}(t)-\mathbf{R}_l(t)\right), \quad t=0,1,\cdots,D-1.
\end{equation}
with $p(t)$ a sequence of $D$ distinct integers comprise between $2$ and $L$.

\section{The Distributed ALS algorithm}\label{sec_DALS}
Recently, a great effort has been devoted to the derivation of distributed estimation algorithms \cite{XiBoKi:07,MaScGi:07,BoDFScVa:08}. Most of them make use of average consensus based estimations. By adopting the average consensus framework, we now derive a distributed version of the alternating least squares (ALS) algorithm, herein referred to as DALS. It is worth mentioning that the formulation of the DALS algorithm proposed here generalizes that of \cite{KiAl:10}. More specifically, herein we assume that each node has a third-order data tensor at its disposal, in contrast to \cite{KiAl:10} where each node is restricted to access a matrix slice of the global tensor. In that case, the developments given in \cite{KiAl:10} can be obtained from the equalities in (\ref{particular_case}) as a particular case.

Given the tensor $\mathbb{X}$, various algorithms can be used for estimating the factor matrices $\mathbf{A}$, $\mathbf{B}$, and $\mathbf{C}$, ALS being the most popular one. By exploiting the unfolded representations of the PARAFAC model (\ref{unfolded_1})-(\ref{unfolded_3}), ALS acts by alternately minimizing the following cost functions $\mathcal{J}_1=\left\|\mathbf{X}_1-\left(\mathbf{B}\odot\mathbf{C}\right)\mathbf{A}^T\right\|_F^2$, $\mathcal{J}_2=\left\|\mathbf{X}_2-\left(\mathbf{C}\odot\mathbf{A}\right)\mathbf{B}^T\right\|_F^2$, and $\mathcal{J}_3=\left\|\mathbf{X}_3-\left(\mathbf{A}\odot\mathbf{B}\right)\mathbf{C}^T \right\|_F^2$.

\subsection{Consensus based estimations of $\mat B$ and $\mat C$}
Let us define the matrix \[\mathbf{Z}=\mathbf{A}\odot\mathbf{C}=\left(
             \begin{array}{c}
               \mathbf{Z}^{(1)}\\
               \vdots \\
              \mathbf{Z}^{(L)}\\
             \end{array}
           \right),
\] with $\mathbf{Z}^{(l)}=\mathbf{A}^{(l)}\odot\mathbf{C}$. From the definition of $\mathbf{X}_2$ we can state that:
\[
\mathbf{X}_2=\pmb{\Pi}\left(\begin{array}{c}
               \mathbf{X}^{(1)}_{2}\\
               \vdots \\
              \mathbf{X}^{(L)}_{2}\\
             \end{array}
           \right)=\pmb{\Pi}_1\mathbf{ZB}^T.
\]
We can also easily show that it exists a permutation matrix $\bar{\pmb{Pi}}$ such that $\pmb{\Pi}_1^T\pmb{\Pi}=\mathbf{I}\otimes\bar{\pmb{\Pi}}$.
As a consequence
\[\mathcal{J}_2=\left\|\mathbf{X}_2-(\mathbf{C}\odot\mathbf{A})\mathbf{B}^T\right\|_F^2
=\left\|\left(\mathbf{I}\otimes\bar{\pmb{\Pi}}\right)\left(\begin{array}{c}
               \mathbf{X}^{(1)}_{2}\\
               \vdots \\
              \mathbf{X}^{(L)}_{2}\\
             \end{array}
           \right)-\mathbf{ZB}^T\right\|_F^2.
\]
Minimizing $\mathcal{J}_2$ as a function of matrix $\mathbf{B}$ yields \[\hat{\mathbf{B}}^T=\left(\mathbf{Z}^H\mathbf{Z}\right)^{-1}\mathbf{Z}^H(\mathbf{I}\otimes\bar{\pmb{\Pi}})\left(\begin{array}{c}
               \mathbf{X}^{(1)}_{2}\\
               \vdots \\
              \mathbf{X}^{(L)}_{2}\\
             \end{array}
           \right).\]
This solution can also be written as:
\begin{equation}
\label{solution_3}
\hat{\mathbf{B}}^T=\left(\frac{1}{L}\sum\limits_{l=1}^{L}\mathbf{Z}^{(l)H}\mathbf{Z}^{(l)}\right)^{-1}\left(\frac{1}{L}\sum\limits_{l=1}^{L}\mathbf{Z}^{(l)H}\bar{\pmb{\Pi}}\mathbf{X}^{(l)}_{2}\right).
\end{equation}
The computation of (\ref{solution_3}) results on averaging local estimates
$\pmb{\Lambda}_l(0)=\mathbf{Z}^{(l)H}\mathbf{Z}^{(l)}$ and $\pmb{\Psi}_{l}(0)=\mathbf{Z}^{(l)H}\bar{\pmb{\Pi}}\mathbf{X}^{(l)}_{2}$. Such an averaging can be achieved using the consensus algorithm. Indeed, we have to run two
average consensus in parallel so that
\[
\pmb{\Lambda}_l(t)\rightarrow\frac{1}{L}\sum\limits_{l=1}^{L}\pmb{\Lambda}_l(0)=\frac{1}{L}\sum\limits_{l=1}^{L}\mathbf{Z}^{(l)H}\mathbf{Z}^{(l)},\]
\[
\pmb{\Psi}_l(t)\rightarrow\frac{1}{L}\sum\limits_{l=1}^{L}\pmb{\Psi}_l(0)=\frac{1}{L}\sum\limits_{l=1}^{L}\mathbf{Z}^{(l)H}\bar{\pmb{\Pi}}\mathbf{X}^{(l)}_{2}.
\]
Therefore, the local estimate of $\mathbf{B}^T$, given by $\hat{\mathbf{B}}^{(l)T}=\pmb{\Lambda}_l^{-1}(t)\pmb{\Psi}_l(t)$ converges towards $\hat{\mathbf{B}}^T$.

The factor matrix $\mathbf{C}$ can be estimated by minimizing the cost function $\mathcal{J}_3$, which can be rewritten as:
\[\mathcal{J}_3=\left\|\mathbf{X}_3-\mathbf{Y}\mathbf{C}^T\right\|_F^2,\] with $\mathbf{Y}=\mathbf{A}\odot\mathbf{B}$. From the property (\ref{prop_khatri_block}) of the Khatri-Rao product, we get:
\[
\mathbf{Y}=\left(
             \begin{array}{c}
               \mathbf{Y}^{(1)}\\
               \vdots \\
              \mathbf{Y}^{(L)}\\
             \end{array}
           \right)=\left(
             \begin{array}{c}
               \mathbf{A}^{(1)}\odot\mathbf{B}\\
               \vdots \\
              \mathbf{A}^{(L)}\odot\mathbf{B}\\
             \end{array}
           \right).
\] By minimizing $\mathcal{J}_3$, given $\mathbf{Y}$, we get \[\hat{\mathbf{C}}^T=\left(\mathbf{Y}^H\mathbf{Y}\right)^{-1}\mathbf{Y}^H\mathbf{X}_3.\]
This solution can also be written as:
\begin{equation}
\label{solution_1}
\hat{\mathbf{C}}^T=\left(\frac{1}{L}\sum\limits_{l=1}^{L}\mathbf{Y}^{(l)H}\mathbf{Y}^{(l)}\right)^{-1}
\left(\frac{1}{L}\sum\limits_{l=1}^{L}\mathbf{Y}^{(l)H}\mathbf{X}^{(l)}_{3}\right).
\end{equation}
As previously, the computation of (\ref{solution_1}) results on averaging local estimates $\pmb{\Gamma}_l(0)=\mathbf{Y}^{(l)H}\mathbf{Y}^{(l)}$ and $\pmb{\Phi}_{l}(0)=\mathbf{Y}^{(l)H}\mathbf{X}^{(l)}_{3}$. Using the consensus algorithm:

\[
\pmb{\Gamma}_l(t)\rightarrow\frac{1}{L}\sum\limits_{l=1}^{L}\pmb{\Gamma}_l(0)=\frac{1}{L}\sum\limits_{l=1}^{L}\mathbf{Y}^{(l)H}\mathbf{Y}^{(l)},\]
\[
\pmb{\Phi}_l(t)\rightarrow\frac{1}{L}\sum\limits_{l=1}^{L}\pmb{\Phi}_l(0)=\frac{1}{L}\sum\limits_{l=1}^{L}\mathbf{Y}^{(l)H}\mathbf{X}^{(l)}_{3}.
\]
Therefore, the local estimates of $\mathbf{C}^T$, defined as $\hat{\mathbf{C}}_l^T=\pmb{\Gamma}_l^{-1}(t)\pmb{\Phi}_l(t)$, $l=1,\cdots,L$, converges towards $\hat{\mathbf{C}}^T$

\subsection{Local estimation of $\mat A^{(l)}$}
Recall that the matrix $\mat A^{(l)}$ is assumed to contain intrinsically local parameters. Therefore there is no need to share these parameters between different nodes. From the definition of the unfolded matrix $\mathbf{X}_1$ and using property (\ref{prop_frobenius2}) of the Frobenius norm, we get:
\[\mathcal{J}_1=\sum\limits_{l=1}^L\left\|\mathbf{X}^{(l)}_{1}-(\mathbf{B}\odot\mathbf{C})\mathbf{A}^{(l)T}\right\|_F^2.\]

As a consequence, the local parameters can be estimated as follows:
\begin{equation}\label{solution_2}
\hat{\mathbf{A}}^{(l)T}=(\hat{\mathbf{B}}\odot\hat{\mathbf{C}})^\dag\mathbf{X}^{(l)}_{1}.
\end{equation}
Notice that the Kruskal condition on the overall tensor guarantees that $\mathbf{B}\odot\mathbf{C}$ is full column rank.

\subsection{Distributed ALS algorithm using average consensus}
The DALS algorithm is constituted by interlacing local ALS steps with consensus iterations. By considering perfect exchanges between nodes during
consensus iterations, the algorithm is summarized below.

Given the sub-tensors $\mathbb{X}^{(l)}$ with unfolded matrices $\mathbf{X}^{(l)}_{i}$, $i=1,2,3$, $l=1,2,\ldots,L$, and the permutation matrix $\bar{\pmb{\Pi}}$:
\begin{enumerate}
\item (\textbf{Initializations}) For $l=1,\cdots,L$, initialize $\hat{\mathbf{A}}^{(l)}(0)$ and $\hat{\mathbf{C}}^{(l)}(0)$ with random values. Set the DALS iteration $i=0$, and select the weights $w_{l,j}$ of the consensus algorithm along with a number $T$ of iterations.
\item\label{D_ALS_loop}(\textbf{Local computations})
For $l=1,\ldots,L$, compute
$\mathbf{Z}^{(l)}(i)=\hat{\mathbf{A}}^{(l)}(i)\odot\hat{\mathbf{C}}^{(l)}(i+1)$,
$\pmb{\Lambda}_l(i,0)=\mathbf{Z}^{(l)H}(i)\mathbf{Z}^{(l)}(i)$, and
$\pmb{\Psi}_l(i,0)=\mathbf{Z}^{(l)H}(i)\bar{\pmb{\Pi}}\mathbf{X}^{(l)}_{2}$.
\item (\textbf{Average consensus}) Run the consensus algorithm for $\pmb{\Lambda}_l$ and $\pmb{\Psi}_l$
\begin{enumerate}
\item For $t=0,1,\cdots,N_c-1$,
\[\pmb{\Lambda}_l(i,t+1)=\pmb{\Lambda}_l(i,t)+\sum\limits_{l'\in\mathcal{N}_l}w_{l,l'}\left(\pmb{\Lambda}_{l'}(i,t)-\pmb{\Lambda}_l(i,t)\right),\]
\[\pmb{\Psi}_l(i,t+1)=\pmb{\Psi}_l(i,t)+\sum\limits_{j\in\mathcal{N}_l}w_{l,j}\left(\pmb{\Psi}_{l'}(i,t)-\pmb{\Psi}_l(i,t)\right).\]
\item Set
$\pmb{\Lambda}_l(i)=\pmb{\Lambda}_l(i,N_c)$ and
$\pmb{\Psi}_l(i)=\pmb{\Psi}_l(i,N_c)$.
\end{enumerate}
\item (\textbf{Local update of the first global matrix}) Compute the local estimates of matrix $\mathbf{B}$
\[
\hat{\mathbf{B}}^{(l)}(i+1)=\pmb{\Lambda}_l^{-1}(i)\pmb{\Psi}_l(i).
\]
\item (\textbf{Local computations}) For $l=1,\ldots,L$, compute
$\mathbf{Y}^{(l)}(i)=\hat{\mathbf{A}}^{(l)}(i)\odot\hat{\mathbf{B}}^{(l)}(i)$,
$\pmb{\Gamma}_l(i,0)=\mathbf{Y}^{(l)H}(i)\mathbf{Y}^{(l)}(i)$, and
$\pmb{\Phi}_l(i,0)=\mathbf{Y}^{(l)H}(i)\mathbf{X}^{(l)}_{3}$.
\item (\textbf{Average consensus}) Run the consensus algorithm for $\pmb{\Gamma}_l$ and $\pmb{\Phi}_l$
\begin{enumerate}
\item For $t=0,1,\cdots,N_c-1$,
\[\pmb{\Gamma}_l(i,t+1)=\pmb{\Gamma}_l(i,t)+\sum\limits_{l'\in\mathcal{N}_l}w_{l,l'}\left(\pmb{\Gamma}_{l'}(i,t)-\pmb{\Gamma}_l(i,t)\right),\]
\[\pmb{\Phi}_l(i,t+1)=\pmb{\Phi}_l(i,t)+\sum\limits_{l'\in\mathcal{N}_l}w_{l,l'}\left(\pmb{\Phi}_{l'}(i,t)-\pmb{\Phi}_l(i,t)\right).\]
\item Set
$\pmb{\Gamma}_l(i)=\pmb{\Gamma}_l(i,N_c)$ and $\pmb{\Phi}_l(i)=\pmb{\Phi}_l(i,N_c)$.
\end{enumerate}
\item (\textbf{Local update of the second global matrix})Compute the local estimates of matrix $\mathbf{C}$
\[
\hat{\mathbf{C}}^{(l)}(i+1)=\pmb{\Gamma}_l^{-1}(i)\pmb{\Phi}_l(i).
\]

\item (\textbf{Local update of the local matrices})Compute the local estimates of $\mat A^{(l)}$
\[
\hat{\mathbf{A}}^{(l)T}(i+1)=\left(\hat{\mathbf{B}}^{(l)}(i)\odot\hat{\mathbf{C}}^{(l)}(i)\right)^\dag\mathbf{X}^{(l)}_{1}.
\]
\item Increment $i$ and return to step \ref{D_ALS_loop} until a convergence criterion is reached.
\end{enumerate}

\section{The Distributed LM algorithm}\label{sec_DLM}
The LM algorithm, originally proposed by Levenberg and Marquardt \cite{Marquardt63}, also known as the ``damped Gauss-Newton method'', has been successfully applied to estimate the parameters of some tensor models. Conceptually, the LM algorithm is mainly based on a modified version of the Gauss-Newton (GN) method, the basic difference being on the insertion of an adjustable regularization parameter (damping factor) into the GN update equations to yield improved convergence properties. The work \cite{ToBr:06} developed an LM algorithm for fitting the PARAFAC model. In \cite{NiDL:08}, an LM algorithm is proposed for estimating the parameters of block tensor models with application to blind multiuser signal separation in code division multiple access (CDMA) systems. Herein, we derive a distributed version of the LM algorithm for in-network computation of the PARAFAC decomposition. Before formulating the distributed LM algorithm for fitting the PARAFAC model, we briefly recall the principle of its conventional (centralized) version.

Define
\begin{equation}\label{eq:vector_p}
\mat p=\left(\begin{array}{c} \textit{vec}(\mat A^T) \\ \textit{vec}(\mat B^T) \\ \textit{vec}(\mat C^T)\end{array}\right)
=\left(\begin{array}{c} \textit{vec}(\mat A^{(1)T}) \\ \vdots \\ \textit{vec}(\mat A^{(L)T})\\ \textit{vec}(\mat B^T) \\ \textit{vec}(\mat C^T)\end{array}\right) \inmat{F}{1}
\end{equation}
with $F=R\big((\sum\limits_{l=1}^LI_l)+J+K\big)$, which concatenates all the model unknowns. The global cost function (\ref{eq:global1}) can then be written as:
\begin{eqnarray}
\mathcal{J}(\mat A,\mat B,\mat
C)&=&\frac{1}{2}\sum\limits_{l=1}^{L}\sum\limits_{i=1}^{I_l}\sum\limits_{j=1}^{J}
\sum\limits_{k=1}^{K}\left|x_{i,j,k}^{(l)}- \sum\limits_{r=1}^{R}
a^{(l)}_{i,r}b_{j,r}c_{k,r}\right|^2\nonumber\\
&=&\frac{1}{2}\sum\limits_{l=1}^{L}\sum\limits_{i=1}^{I_l}\sum\limits_{j=1}^{J}\sum\limits_{k=1}^{K}\left|
r^{(l)}_{i,j,k}(\mat A,\mat B,\mat
C)\right|^2\nonumber\\
&=& \frac{1}{2}\sum\limits_{l=1}^{L}\left\| \mat r_l(\mat p)\right\|^2\label{eq:cost_residuals}
\end{eqnarray}
where 
$\mat r_l(\mat p) \inmat{I_lJK}{1}$ is the vector of residuals associated with node $l$, the $m$-th element of which, $[\mat r_l(\mat p)]_m$, is given by $[\mat r_l(\mat p)]_m= r^{(l)}_{i,j,k}(\mat A,\mat B,\mat C)$ with $m=(k-1)JI+ (j-1)I +i$.
By making use of the unfolded representations of the PARAFAC model given in (\ref{eq:unfold_X1l})-(\ref{eq:unfold_X3l}), we have:
\beq
\mat r_l(\mat p)&=& \textit{vec}(\mat X^{(l)}_1) - \Big[\mat I_{I_l} \otimes (\mathbf{B}\odot\mathbf{C})\Big]\textit{vec}(\mathbf{A}^{(l)T}) \,\, \in\mathbb{C}^{I_lJK\times 1},\nonumber\\
&=& \bm{\Pi}_2\Big(\textit{vec}(\mat X^{(l)}_2) -  \Big[\mat I_{J} \otimes (\mathbf{C}\odot\mathbf{A}^{(l)})\Big]\textit{vec}(\mathbf{B}^T)\Big) \,\, \in\mathbb{C}^{I_lJK\times 1},\nonumber\\
&=&\bm{\Pi}_3\Big(\textit{vec}(\mat X^{(l)}_3) -  \Big[\mat I_{K} \otimes (\mathbf{A}^{(l)}\odot\mathbf{B})\Big]\textit{vec}(\mathbf{C}^T) \Big)\,\, \in\mathbb{C}^{I_lJK\times 1}\nonumber,
\eeq
where $\bm{\Pi}_2 \inmat{I_lJK}{JKI_l}$ and $\bm{\Pi}_3 \inmat{I_lJK}{KI_lJ}$ are permutation matrices.


The principle of the GN algorithm is to approximate the vector of residuals in the neighborhood of $\mat p$ by a Taylor expansion truncated after the first-order term. Therefore, for $\|\Delta\mat{p}\|$ sufficiently small, we have:
\begin{equation}\label{eq:lin_approx}
\mat r(\mat p+\Delta \mat p)\cong \mat r(\mat p) + \mat J\Delta \mat p,
\end{equation}
where 
$\mat J \inmat{(\sum\limits_{l=1}^LI_l)JK}{F}$ is the Jacobian matrix.
Using the linear approximation (\ref{eq:lin_approx}), we can rewrite
(\ref{eq:cost_residuals}) in terms of $\Delta \mat p$, as:
\begin{equation}\label{eq:cost_Deltap}
\tilde{\mathcal{J}}(\Delta \mat p)=\left\|\mat r(\mat p) + \mat J\Delta \mat p\right\|^2.
\end{equation}
The correction term $\Delta \mat p$ is then computed as the solution of the linear least squares problem:
\[\underset{\Delta \mat p}{\emph{min}}\left\|\mat r(\mat p) + \mat
J\Delta \mat p\right\|^2,
\]
which is given by the system of normal equations:
\begin{equation}\label{eq:normal_eqs}
(\mat J^{H}\mat J)\Delta \mat p=-\mat J\mat r(\mat p),
\end{equation}
where 
$\mat J^H\mat J$ is an approximation to the Hessian matrix \cite{MadsNT04}. The idea of the LM algorithm is to add a regularization parameter (damping factor)$\footnote{We refer the interested reader to \cite{MadsNT04} for further details and discussions on the damping factor $\lambda$.}$ to (\ref{eq:normal_eqs}), thus yielding:
\begin{equation}\label{eq:normal_eqs_LM}
(\mat J^{H}\mat J+\lambda\mat I_F)\Delta \mat p=-\mat J^{H}\mat r(\mat p). 
\end{equation}

\subsection{Partioned estimation of local and global parameters}
From the partitioned structure of $\mat p$ in (\ref{eq:vector_p}), we can write the Jacobian $\mat J$ as a concatenation of matrix
blocks in the following way:
\begin{equation}\label{eq:Jac_global}
\mat J=\left(\begin{array}{ccccc} \frac{\partial \mat r_1(\mat p)}{\partial \textit{vec}(\mat A^{(1)T}) } & \cdots &  \frac{\partial \mat r_1(\mat p)}{\partial \textit{vec}(\mat A^{(L)T}) } & \frac{\partial \mat r_1(\mat p)}{\partial \textit{vec}(\mat B^T)}  & \frac{\partial \mat r_1(\mat p)}{\partial \textit{vec}(\mat C^T)}\\
\vdots & \vdots & \cdots & \vdots & \vdots\\
\frac{\partial \mat r_L(\mat p)}{\partial \textit{vec}(\mat A^{(1)T}) } & \cdots &  \frac{\partial \mat r_L(\mat p)}{\partial \textit{vec}(\mat A^{(L)T})}& \frac{\partial \mat r_1(\mat p)}{\partial \textit{vec}(\mat B^T)}  & \frac{\partial \mat r_1(\mat p)}{\partial \textit{vec}(\mat C^T)}
\end{array}\right)
\end{equation}
where
\begin{eqnarray}
&&\frac{\partial \mat r_l(\mat p)}{\partial \textit{vec}(\mat A^{(l')T}) }=- \Big[\mat I_{I_l} \otimes (\mathbf{B}\odot\mathbf{C})\Big]\delta_{l,l'}\inmat{I_lJK}{RI_l},\label{eq:jac1}\\
&&\frac{\partial \mat r_l(\mat p)}{\partial \textit{vec}(\mat B^T)}= -  \bm{\Pi}_2\Big[\mat I_{J} \otimes (\mathbf{C}\odot\mathbf{A}^{(l)})\Big]\inmat{I_lJK}{RJ},\label{eq:jac2}\\
&&\frac{\partial \mat r_l(\mat p)}{\partial \textit{vec}(\mat C^T)}= -  \bm{\Pi}_3\Big[\mat I_{K} \otimes (\mathbf{A}^{(l)}\odot\mathbf{B})\Big]\inmat{I_lJK}{RK}.\label{eq:jac3}
\end{eqnarray}
Let us partition the parameter vector $\mat p$ into $L+1$ sub-vectors, such that $\overline{\mat p}= [\mat a^{T}_1,\ldots, \mat a^{T}_L,\overline{\mat p}^T]^T$, where
$\mat a_l=\textit{vec}(\mat A^{(l)T})$, $l=1,\ldots, L$, and $\overline{\mat p}=[\textit{vec}(\mat B^T)^T \textit{vec}(\mat C^T)^T]^T \inmat{R(J+K)}{1}$.
From this partitioning, the Jacobian matrix (\ref{eq:Jac_global}) can be rewritten as:
\[
\mat J(\mat p)=\left(\begin{array}{cccc}\mat J(\mat a_1) &  & &
\mat J_1(\overline{\mat p})\\   & \ddots & & \vdots\\ & & \mat J(\mat a_L)  &
\mat J_{L}(\overline{\mat p})
\end{array}\right),
\]
where
\be
\mat J(\mat a_l) =\frac{\partial \mat r_l(\mat p)}{\partial \textit{vec}(\mat A^{(l)T})}\inmat{I_lJK}{RI_l},\quad \mat J_{l}(\overline{\mat p})=\left[\frac{\partial \mat r_l(\mat p)}{\partial \textit{vec}(\mat B^T)}\,\, , \,\,\frac{\partial \mat r_l(\mat p)}{\partial \textit{vec}(\mat C^T)}\right]\inmat{I_lJK}{R(J+K)}.\label{eq:jac_partitioned}
\ee
Therefore, the approximation of the Hessian matrix (left hand side of (\ref{eq:normal_eqs})) is given by:
\[
\mat J^H(\mat p)\mat J(\mat p)=\left(\begin{array}{cccc}\mat J^H(\mat
a_1)\mat J(\mat a_1)& &  & \mat J^H(\mat a_1)\mat J_1(\overline{\mat p})\\  &
\ddots & & \vdots\\  & & \mat J^H(\mat a_L)\mat J(\mat a_L) &
\mat J^H(\mat a_L)\mat J_{L}(\overline{\mat p})\\
\mat J^H_{1}(\overline{\mat p})\mat J(\mat a_1) & \cdots &
\mat J^H_{L}(\overline{\mat p})\mat J(\mat a_L) &
\sum\limits_{l=1}^{L}\mat J^H_{l}(\overline{\mat
p})\mat J_{l}(\overline{\mat p})
\end{array}\right),
\]
or, equivalently,
\begin{equation}\label{eq:JHJ_equivalent}
\mat J^H(\mat p)\mat J(\mat p)=\left(\begin{array}{cccc} \mat H(\mat a_1)& &
& \mat Q^H_1\\  & \ddots & & \vdots\\  & & \mat H(\mat a_L) &
\mat Q^H_{L}\\
\mat Q_1 & \cdots & \mat Q_{L} &
\sum\limits_{l=1}^{L}\mat H_l(\overline{\mat p})
\end{array}\right),
\end{equation}
where
\beq
&&\mat H(\mat a_l)=\mat J^H(\mat a_l)\mat J(\mat a_l) \inmat{RI_l}{RI_l},\quad \mat H_l(\overline{\mat p})=\mat J^H_{l}(\overline{\mat p})\mat J_{l}(\overline{\mat p}) \inmat{R(J+K)}{R(J+K)},\nonumber\\
&&\mat B_l=\mat J^H_{l}(\overline{\mat p})\mat J(\mat a_l) \inmat{R(J+K)}{RI_l}.\nonumber
\eeq
The expressions above are detailed in subsection \ref{subsec_detail}.

For the right-hand side of (\ref{eq:normal_eqs}), we obtain:
\begin{eqnarray}\label{eq:JHr_equivalent}
\mat J^H(\mat p)\mat r(\mat p)=\left(\begin{array}{cc} \mat J^H(\mat
a_1)\mat r_1(\mat p) \\ \vdots \\\mat J^H(\mat
a_L)\mat r_L(\mat p)\\
\sum\limits_{l=1}^{L}\mat J^H_l(\overline{\mat p})\mat r_l(\mat p)\end{array}\right). 
\end{eqnarray}

Using (\ref{eq:JHJ_equivalent}) and (\ref{eq:JHr_equivalent}), the normal equations (\ref{eq:normal_eqs}) for the update of the parameter vector $\mat p$ 
can be reformulated as:
\[\left(\begin{array}{cccc} \mat H(\mat a_1)& &
& \mat Q^H_1\\  & \ddots & & \vdots\\  & & \mat H(\mat a_L) &
\mat Q^H_{L}\\
\mat Q_1 & \cdots & \mat Q_{L} &
\sum\limits_{l=1}^{L}\mat H_l(\overline{\mat p})
\end{array}\right)\left(\begin{array}{c} \Delta\mat a_1 \\ \vdots \\ \Delta\mat a_L \\ \Delta \overline{\mat p}\end{array}\right) =-\left(\begin{array}{cc} \mat J^H(\mat
a_1)\mat r_1(\mat p) \\ \vdots \\\mat J^H(\mat
a_L)\mat r_L(\mat p)\\
\sum\limits_{l=1}^{L}\mat J^H_l(\overline{\mat p})\mat r_l(\mat p)\end{array}\right),
\]
which is equivalent to
\begin{equation}\label{eq:normal_eqs_partitioned}
\left(\begin{array}{cc} \mat H(\mat a_1)\Delta\mat a_1 +
\mat Q^H_1\Delta\overline{\mat p}\\ \vdots\\
\mat H(\mat a_L)\Delta\mat a_L + \mat
Q^H_{L}\Delta\overline{\mat p}\\
\sum\limits_{l=1}^{L}\left(\mat Q_l\Delta\mat a_l +
\mat H_l(\overline{\mat p})\Delta\overline{\mat
p}\right)\end{array}\right)=-\left(\begin{array}{cc} \mat J^H(\mat a_1)\mat
r_1(\mat p) \\ \vdots \\\mat J^H(\mat a_L)\mat r_L(\mat
p)\\ \sum\limits_{l=1}^{L}\mat J^H_l(\overline{\mat p})\mat r_l(\mat p)\end{array}\right).
\end{equation}
From (\ref{eq:normal_eqs_partitioned}), we get:
\beq
&&\mat H(\mat a_l)\Delta\mat a_l + \mat Q^H_l\Delta\overline{\mat p}= - \mat J^H(\mat a_l)\mat
r_l(\mat p),\quad l=1,\ldots, L,\label{eq:Deltaak}\\
&&\sum\limits_{l=1}^{L}\left(\mat Q_l\Delta\mat a_l + \mat H_l(\overline{\mat p})\Delta\overline{\mat p}\right)= - \sum\limits_{l=1}^{L}\mat J^H_l(\overline{\mat p})\mat r_l(\mat p)\label{eq:Deltap}.
\eeq
In (\ref{eq:Deltaak}), the correction term $\Delta\mat a_l$ can be written as:
\begin{equation}\label{eq:Deltaak2}
\Delta\mat a_l =- \mat H^{-1}(\mat a_l)\left(\mat Q^H_l\Delta\overline{\mat p} +  \mat J^H(\mat a_l)\mat r_l(\mat p)\right).
\end{equation}
Then, substituting (\ref{eq:Deltaak2}) into (\ref{eq:Deltap}) yields:
\beq
&&\sum\limits_{l=1}^{L}\left( -\mat Q_l\mat H^{-1}(\mat a_l)\mat Q^H_l\Delta\overline{\mat p} - \mat Q_l\mat H^{-1}(\mat a_l)\mat J^H(\mat a_l)\mat r_l(\mat p) + \mat H_l(\overline{\mat p})\Delta\overline{\mat p}\right)\nonumber\\
&&= - \sum\limits_{l=1}^{L}\mat J^H_l(\overline{\mat p})\mat r_l(\mat p).\nonumber\label{eq:Deltap2}
\eeq
that can be compactly written as
\be
\left(\sum\limits_{l=1}^{L}\bm{\Theta}_l\right)\Delta\overline{\mat p}= - \sum\limits_{l=1}^{L} \bm{\xi}_l\label{eq:Deltap4},
\ee
where
\beq
&&\bm{\Theta}_l=\mat H_l(\overline{\mat p}) - \mat Q_l\mat H^{-1}(\mat a_l)\mat Q^H_l,\\
&& \bm{\xi}_l= -\left(\mat J^H_l(\overline{\mat p}) - \mat Q_l\mat H^{-1}(\mat a_l)\mat J^H(\mat a_l)\right)\mat r_l(\mat p),
\eeq
represent the equivalent Hessian matrix and gradient vector calculated at the node $l$, respectively. Finally, by adding the regularization parameter $\lambda$ to (\ref{eq:Deltap4}), we get the following solution for the step $\Delta\overline{\mat p}$:
\be\label{eq:solution_LM}
\Delta\overline{\mat p}= \left[\sum\limits_{l=1}^{L}\bm{\Theta}_l + \bar{\lambda}\mat I\right]^{-1}\left(\sum\limits_{l=1}^{L} \bm{\xi}_l\right).
\ee

\subsubsection{Consensus-based estimation of global parameters}
From (\ref{eq:solution_LM}), we can note that the computation of $\Delta \overline{\mat p}$ supposes the knowledge of the set of Hessian matrices $\{\bm{\Theta}_1,\ldots, \bm{\Theta}_L\}$ and gradient vectors $\{\bm{\xi}_1,\ldots, \bm{\xi}_L\}$  at each node. Otherwise stated, before computing $\Delta \overline{\mat p}$, the nodes must exchange their estimated Hessian matrices and gradient vectors. In fact, just a in-network computation of the average of local Hessian matrices and gradient vectors is required. Indeed, we can rewrite (\ref{eq:solution_LM}) as:

\be\label{eq:solution_LM_2}
\Delta\overline{\mat p}= \left[\bm{\Theta} + \left(\frac{\lambda}{L}\right)\mat I\right]^{-1}\bm{\xi}
\ee
with $\bm{\Theta}=\frac{1}{L}\sum\limits_{l=1}^{L}\bm{\Theta}_l$, $\bm{\xi}=\frac{1}{L}\bm{\xi}_l$, and $\frac{\lambda}{L}=\bar{\lambda}$.

As for the DALS algorithm, $\bm{\Theta}$ and $\bm{\xi}$ can be computed using an average consensus algorithm:
\begin{eqnarray}\label{eq:averageHk}
&&\bm{\Theta}_l(t+1)= \bm{\Theta}_l(t)+ \sum\limits_{l'\in \matc N_l}w_{l,l'}\left(\bm{\Theta}_{l'}(t)-\bm{\Theta}_l(t)\right),\\
&&\bm{\xi}_l(t+1)= \bm{\xi}_l(t)+ \sum\limits_{l'\in \matc N_l}w_{l,l'}\left(\bm{\xi}_{l'}(t)-\bm{\xi}_l(t)\right).
\end{eqnarray}

Consequently, the step $\Delta\overline{\mat p}_l$ computed at each node $l$:

\be\label{eq:solution_LM_l}
\Delta\overline{\mat p}_l= \left[\bm{\Theta}_l(t) + \left(\frac{\lambda}{L}\right)\mat I\right]^{-1}\bm{\xi}_l(t)
\ee
is asymptotically equal to $\Delta\overline{\mat p}$ given in (\ref{eq:solution_LM}).
\subsubsection{Estimation of local parameters}
Since the parameter vector $\mat a_l$ is strictly local, no exchange of information is required to compute of $\Delta \mat a_l$ at the $l$-th node. Therefore, from (\ref{eq:Deltaak2}), the local parameter vector can be estimated from the following steps:
\beq\label{eq:Deltaak2_v2}
&&\Delta\mat a_l =- \mat H^{-1}(\mat a_l)\left(\mat Q^H_l\Delta\overline{\mat p} +  \mat J^H(\mat a_l)\mat r_l(\mat p)\right),\\
&&\mat a_l(i+1) = \mat a_l(i) + \Delta{\mat a_l}.
\eeq

\subsubsection{Calculations of $\mat H(\mat a_l)$, $\mat H_l(\overline{\mat p})$ and $\mat B_l$}\label{subsec_detail} The expressions of $\mat H(\mat a_l)$, $\mat H_l(\overline{\mat p})$ and $\mat B_l$ involved in the computation of $\bm \Theta_l$ and $\bm \xi_l$ are detailed in this subsection.
From the expressions given in (\ref{eq:jac1})-(\ref{eq:jac3}) and the definitions in (\ref{eq:jac_partitioned}), we get:
\beq
\mat H(\mat a_l)= \mat H^{(aa)}_l&=& - \mat I_{I_l} \otimes (\mat B\odot \mat C)^H(\mat B\odot \mat C)\nonumber\\
&=& - \mat I_{I_l} \otimes(\mat B^H\mat B \star \mat C^H\mat C),\label{eq:defHk}
\eeq
\beq
\mat B_l= \left(\begin{array}{c}\Big[\mat I_{J} \otimes (\mat C \odot \mathbf{A}^{(l)})\Big]^H\bm{\Pi}^T_2\Big[\mat I_{I_l} \otimes (\mathbf{B}\odot\mathbf{C})\Big]\\ \Big[\mat I_{K} \otimes (\mathbf{A}^{(l)} \odot \mat B)\Big]^H\bm{\Pi}^T_3\Big[\mat I_{I_l} \otimes (\mathbf{B}\odot\mathbf{C})\Big]\end{array}\right),
\eeq
and
\be
\mat H_l(\overline{\mat p})=\left(\begin{array}{cc}
 \mat H^{(bb)}_l&\mat H^{(bc)H}_l\\
\mat H^{(bc)}_l& \mat H^{(cc)}_l\end{array}\right),\label{eq:Hbar}
\ee
with
\beq
\mat H^{(bb)}_l&=& \Big[\mat I_{J} \otimes (\mat C \odot \mathbf{A}^{(l)})\Big]^H\bm{\Pi}^T_2\bm{\Pi}_2\Big[\mat I_{J} \otimes (\mat C \odot \mathbf{A}^{(l)})\Big]\nonumber\\
&=& \mat I_{J} \otimes (\mathbf{C}^H\mat C \star \mathbf{A}^{(l)H}\mathbf{A}^{(l)} )\label{eq:Hbb}\\
\mat H^{(bc)}_l&=& \bm{\Pi}_3\Big[\mat I_{K} \otimes (\mathbf{A}^{(l)}\odot\mathbf{B})\Big]^H\Big[\mat I_{J} \otimes (\mathbf{C}\odot\mathbf{A}^{(l)})\Big]\bm{\Pi}_2\label{eq:Hbc}\\
\mat H^{(cc)}_l&=& \Big[\mat I_{K} \otimes (\mathbf{A}^{(l)} \odot \mat B)\Big]^H\bm{\Pi}^T_3\bm{\Pi}_3\Big[\mat I_{K} \otimes (\mathbf{A}^{(l)} \odot \mat B)\Big]\nonumber\\
&=& \mat I_{K} \otimes (\mathbf{A}^{(l)H}\mathbf{A}^{(l)} \star  \mathbf{B}^H\mat B).\label{eq:Hcc}
\eeq
Therefore, the Jacobian matrix at each node $l$ can be constructed from the matrices $\mat H^{(aa)}_l$, $\mat H^{(bb)}_l$, $\mat H^{(bc)}_l$, and $\mat H^{(cc)}_l$ calculated using (\ref{eq:defHk}), (\ref{eq:Hbb}), (\ref{eq:Hbc}), and (\ref{eq:Hcc}), respectively.

\subsection{Distributed LM algorithm using average consensus}
The distributed LM algorithm operates as follows:

%
\begin{enumerate}

\item (\textbf{Initializations}) Set $i=0$. Each node initializes randomly $\hat{\mat A}^{(l)}(i=0)$, $\hat{\mat B}^{(l)}(i=0)$, $\hat{\mat C}^{(l)}(i=0)$; From $\hat{\mat B}^{(l)}(i=0)$ and $\hat{\mat C}^{(l)}(i=0)$, build $\hat{\overline{\mat p}}(i=0)$. 
    Choose an initial value for the damping parameter $\lambda$;

\item $i\leftarrow i+1$;

\item (\textbf{Local computations}) Each node computes:

\begin{itemize}

\item Compute $\mat J(\hat{\mat a}_l)= \mat I_{I_l} \otimes (\hat{\mathbf{B}}^{(l)}(i)\odot\hat{\mathbf{C}}^{(l)}(i))$;

\item Compute $\mat J_l(\hat{\overline{\mat p}})=- \left(\bm{\Pi}_2\Big[\mat I_{J} \otimes (\hat{\mathbf{C}}(i)\odot\hat{\mathbf{A}}^{(l)})(i)\Big],\,\, \bm{\Pi}_3\Big[\mat I_{K} \otimes (\hat{\mathbf{A}}^{(l)}(i)\odot\hat{\mathbf{B}}(i))\Big]\right)$;

\item Compute $\mat H(\hat{\mat a}_l)= \mat J^H(\hat{\mat a}_l)\mat J(\hat{\mat a}_l)$, $\mat H_l(\hat{\overline{\mat p}})=\mat J^H_l(\hat{\overline{\mat p}})\mat J_l(\hat{\overline{\mat p}})$ and $\mat Q_l=\mat J^H_l(\hat{\overline{\mat p}})\mat J(\hat{\mat a}_l)$ from expressions (\ref{eq:defHk})-(\ref{eq:Hcc});

\item Compute $\bm{\Theta}_l= \mat H_l(\hat{\overline{\mat p}}) - \mat Q_l\mat H^{-1}(\hat{\mat a}_l)\mat Q^H_l$;

\item Compute $\bm{\xi}_l= -\left(\mat J^H_l(\hat{\overline{\mat p}}) - \mat Q_l\mat H^{-1}(\hat{\mat a}_l)\mat J^H(\hat{\mat a}_l)\right)\mat r_l(\hat{\mat p})$;

\end{itemize}

\item (\textbf{Average consensus}) Initialize $\bm{\Theta}_l(0)=\bm{\Theta}_l$ and $\bm{\xi}_l(0)=\bm{\xi}_l$.\\
For $t=0,1, \ldots, N_c-1$
\begin{eqnarray}\label{eq:averageHk}
&&\hspace{-1.5cm}\bm{\Theta}_l(t+1)= \bm{\Theta}_l(t)+
\sum\limits_{l'\in \matc N_l}w_{l,l'}\left(\bm{\Theta}_{l'}(t)-\bm{\Theta}_l(t)\right)\nonumber\\
&&\hspace{-1.5cm}\bm{\xi}_l(t+1)= \bm{\xi}_l(t)+
\sum\limits_{l'\in \matc N_l}w_{l,l'}\left(\bm{\xi}_{l'}(t)-\bm{\xi}_l(t)\right)\nonumber
\end{eqnarray}

\item (\textbf{Parameter updates}) After $N_c$ iterations, each node uses $\bm{\Theta}_l(N_c)$ and $\bm{\xi}_l(N_c)$ to compute:

\beq
&&\Delta{\hat{\overline{\mat p}}}_l= (\bm{\Theta}_l(N_c) +\lambda\mat I)^{-1}\bm{\xi}_l(N_c)\nonumber\\
&&\hat{\overline{\mat p}}_l(i+1) = \hat{\overline{\mat p}}_l(i) + \Delta{\hat{\overline{\mat p}}}_l\nonumber\\
&&\Delta\hat{\mat a}_l =- \mat H^{-1}(\hat{\mat a}_l)\left(\mat Q^H_l\Delta\hat{\overline{\mat p}}_l +  \mat J^H(\hat{\mat a}_l)\mat r_l(\hat{\mat p})\right)\nonumber\\
&&\hat{\mat a}_l(i+1) = \hat{\mat a}_l(i) + \Delta{\hat{\mat a}_l};\nonumber
\eeq

\item Build  the local estimate of $\hat{\mat A}^{(l)}(i+1)$ from $\hat{\mat a}_l(i+1)$ and those of $\hat{\mat B}^{(l)}(i+1)$ and $\hat{\mat C}^{(l)}(i+1)$ from the vector $\hat{\overline{\mat p}}_l(i+1)$.


\item Repeat Steps 2 to 7 until convergence is achieved.

\end{enumerate}

\section{Simulation Results}\label{sec:simres}
In this section, we present some results obtained by simulating a network of $L=9$ nodes whose communication links are modeled with a connected graph assumed to be a cycle or a Payley graph (see \cite{GoRo:01}). Each node has at its disposal the data samples $x_{i,j,k}^{(l)}=\sum\limits_{r=1}^Ra_{i,r}^{(l)}b_{j,r}c_{k,r}$, with $R=4$, $i=1$, $j=1,\cdots,4$, and $k=1,\cdots,10$, meaning that sub-tensors at each node are simply matrices. This scenario is similar to that considered in \cite{KiAl:10} for collaborative blind symbol detection in CDMA networks. The data $c_{k,r}$ and $b_{j,r}$ were random binary data with values $\pm 1$ while $a_{i,r}$ were generated from a unit normal distribution. In what follows, we evaluate the performance of the proposed methods in terms of NMSE (Normalized Mean Square Error) defined as:
\[
NMSE=\frac{1}{L}\sum\limits_{l=1}^L\frac{\left\|\tilde{\mathbf{X}}_3^{(l)}-\left(\hat{\mathbf{A}}^{(l)}\odot\hat{\mathbf{B}}\right)\hat{\mathbf{C}}^T\right\|_F^2}{\left\|\mathbf{X}_3^{(l)}\right|_F^2}
\]
or equivalently as
\[
NMSE=\frac{1}{L}\sum\limits_{l=1}^L\frac{\left\|\tilde{\mathbf{X}}_3^{(l)}-\hat{\mathbf{B}}diag(\hat{\mathbf{A}}^{(l)})\hat{\mathbf{C}}^T\right\|_F^2}{\left\|\mathbf{X}_3^{(l)}\right|_F^2}
\] since each $\mathbf{A}^{(l)}$ is a vector. Here $\tilde{\mathbf{X}}_3^{(l)}$ stands for possibly noisy observations of $\mathbf{X}_3^{(l)}$ at node $l$.
The results presented herein are median values over $100$ independent Monte-Carlo runs. The iterations of the evaluated algorithms (DALS, DLM, and their centralized counterparts) were stopped after $100$ iterations. For DLM, the damping parameter is set equal to $10^{-3}$.
\par We first consider the ideal case where both observations at each node and data exchange during consensus iterations are noiseless. The average consensus protocol used here is the optimal constant edge policy as proposed in \cite{XiBo:04}. The consensus matrix is computed as $\mathbf{W}=\mathbf{I}-\frac{2}{\lambda_2+\lambda_L}\mathbf{L}$, where $\lambda_2$ and $\lambda_L$ are respectively the second smallest and the largest eigenvalues of the graph Laplacian matrix $\mathbf{L}$. Such a policy ensures fast convergence of the average consensus protocol. Since this standard consensus protocol guarantees only an asymptotic convergence, we evaluate the effect of imperfect average consensus by considering a finite number $N_c$ of consensus iterations.

%
%
%

\begin{figure}[!t]
\begin{center}
\includegraphics[scale=0.8]{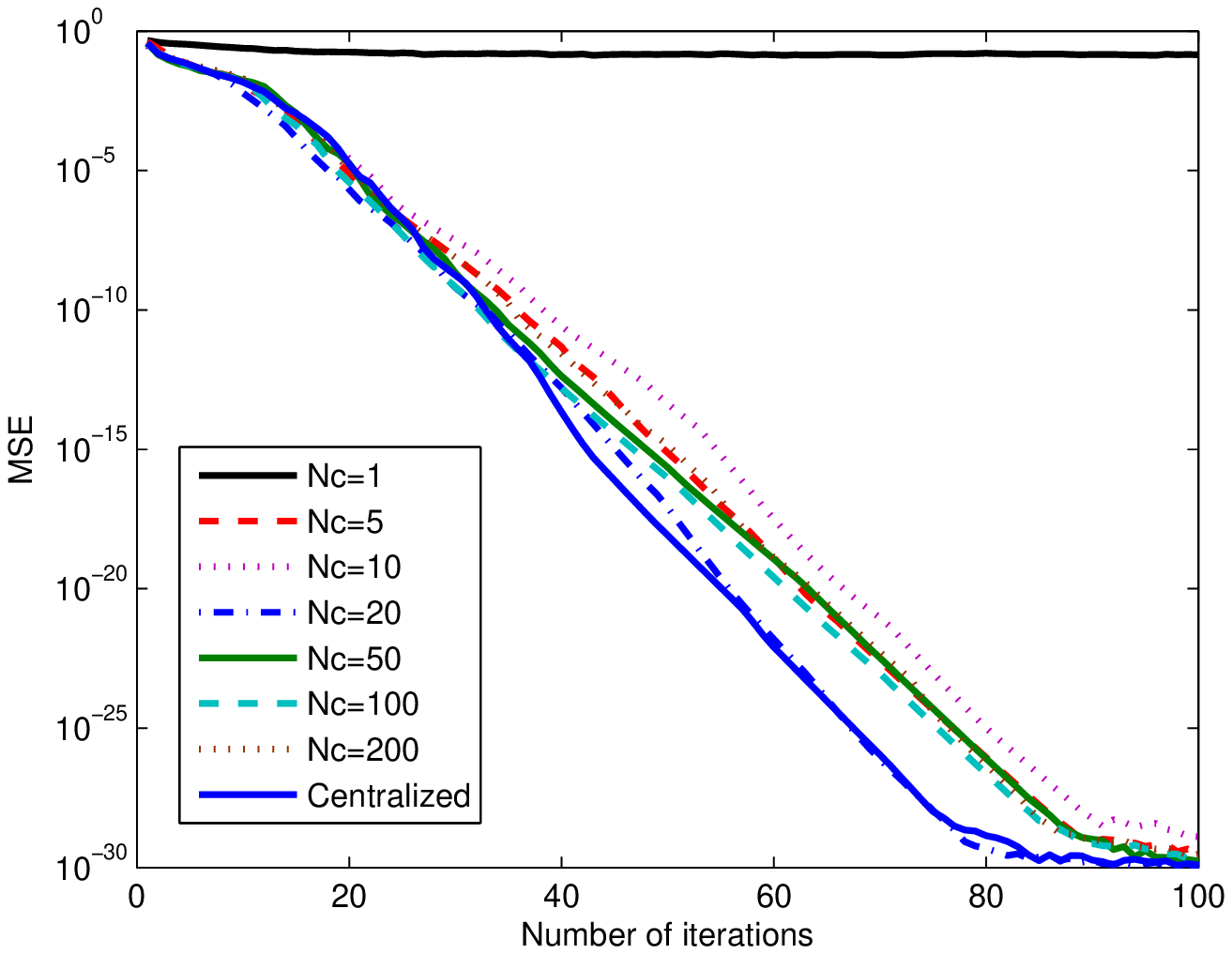}
\end{center}
\caption{Mean square for the DALS algorithm in the case of a cycle graph.}
\label{median_DALS_cycle_fig}
\end{figure}

\begin{figure}[!t]
\begin{center}
\includegraphics[scale=0.8]{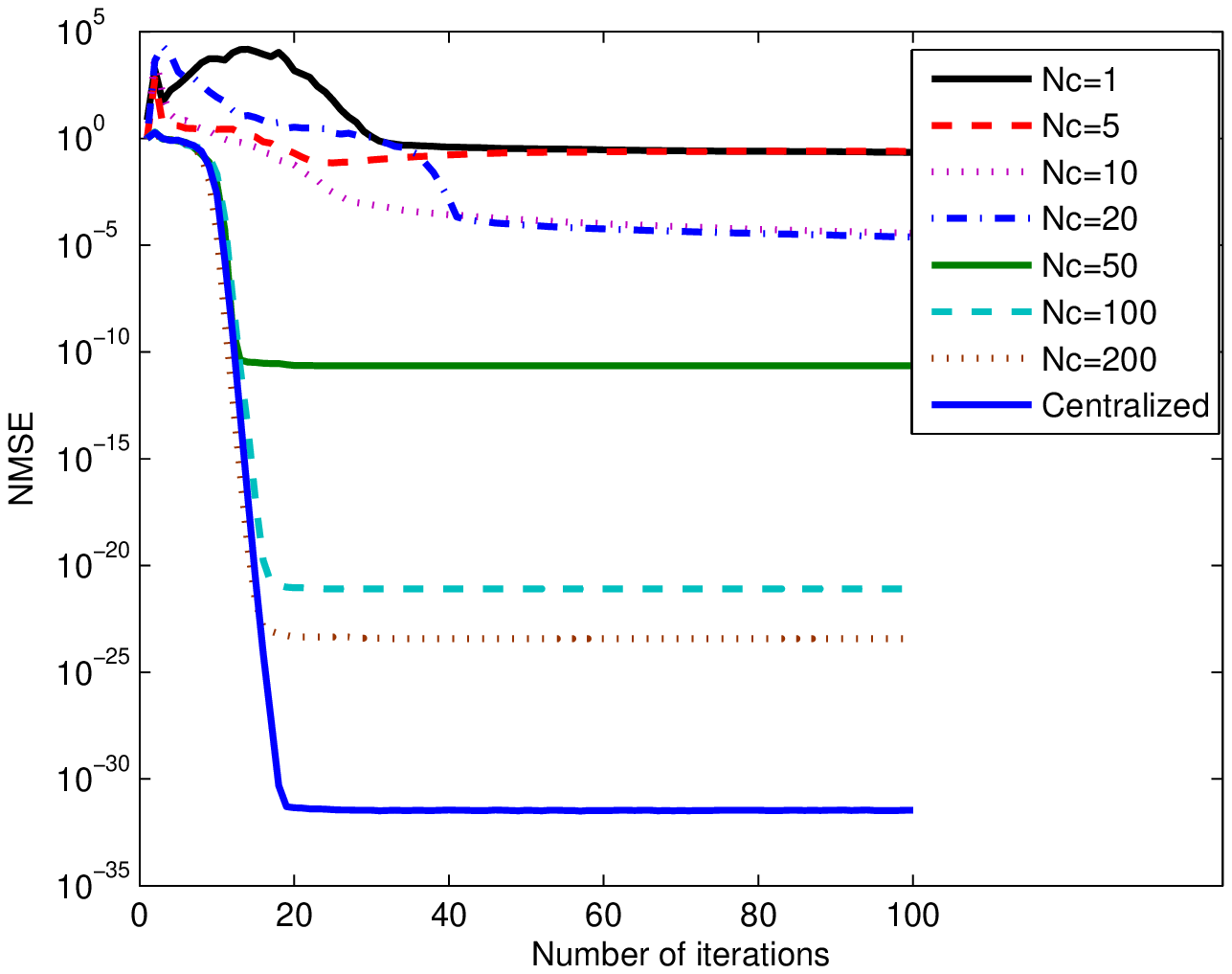}
\end{center}
\caption{Mean square for the DLM algorithm in the case of a cycle graph.}
\label{median_DLM_cycle fig}
\end{figure}

\begin{figure}[!t]
\begin{center}
\includegraphics[scale=0.8]{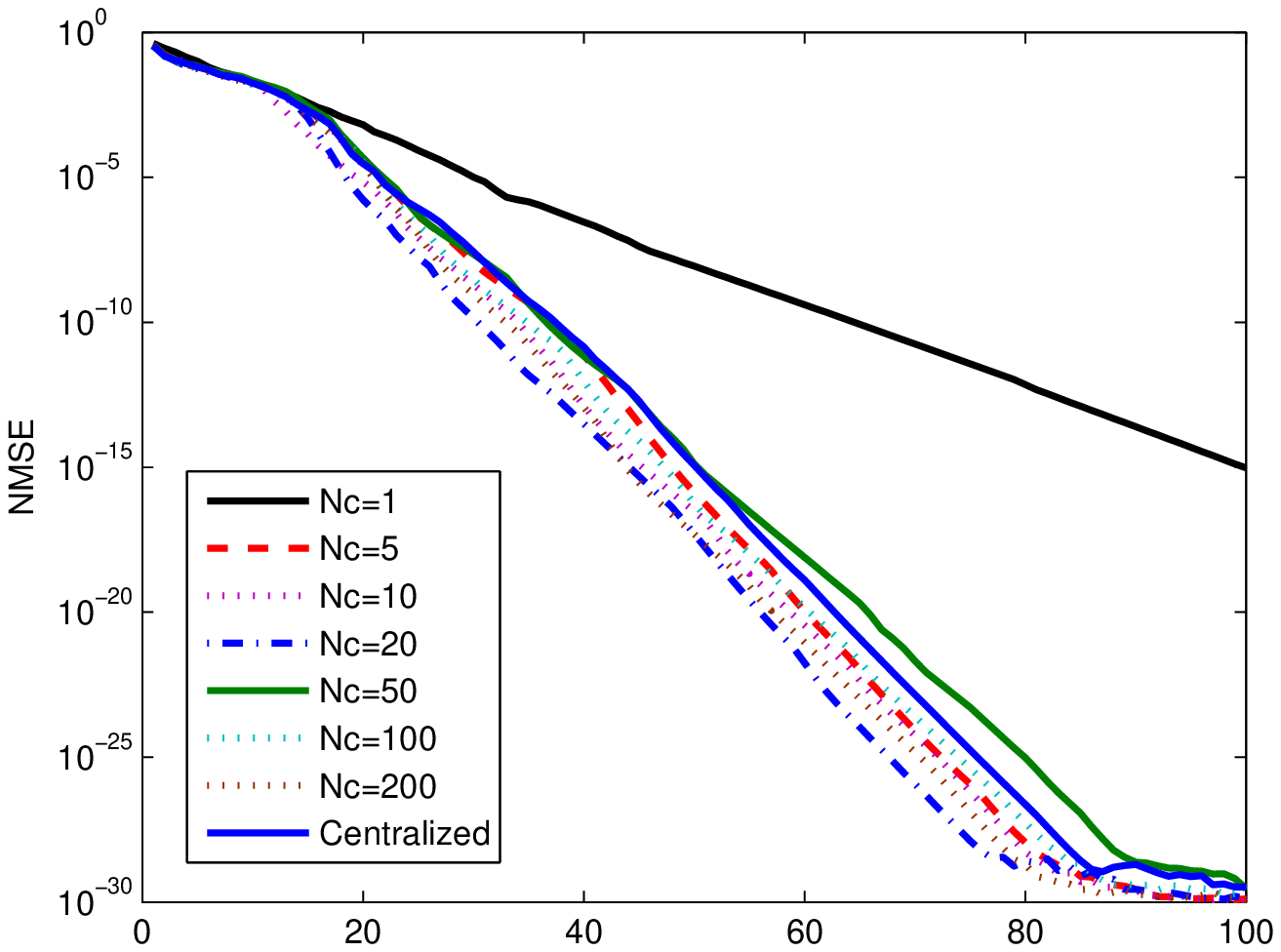}
\end{center}
\caption{Mean square for the DALS algorithm in the case of a Payley graph.}
\label{median_DALS_payley fig}
\end{figure}

\begin{figure}[!t]
\begin{center}
\includegraphics[scale=0.8]{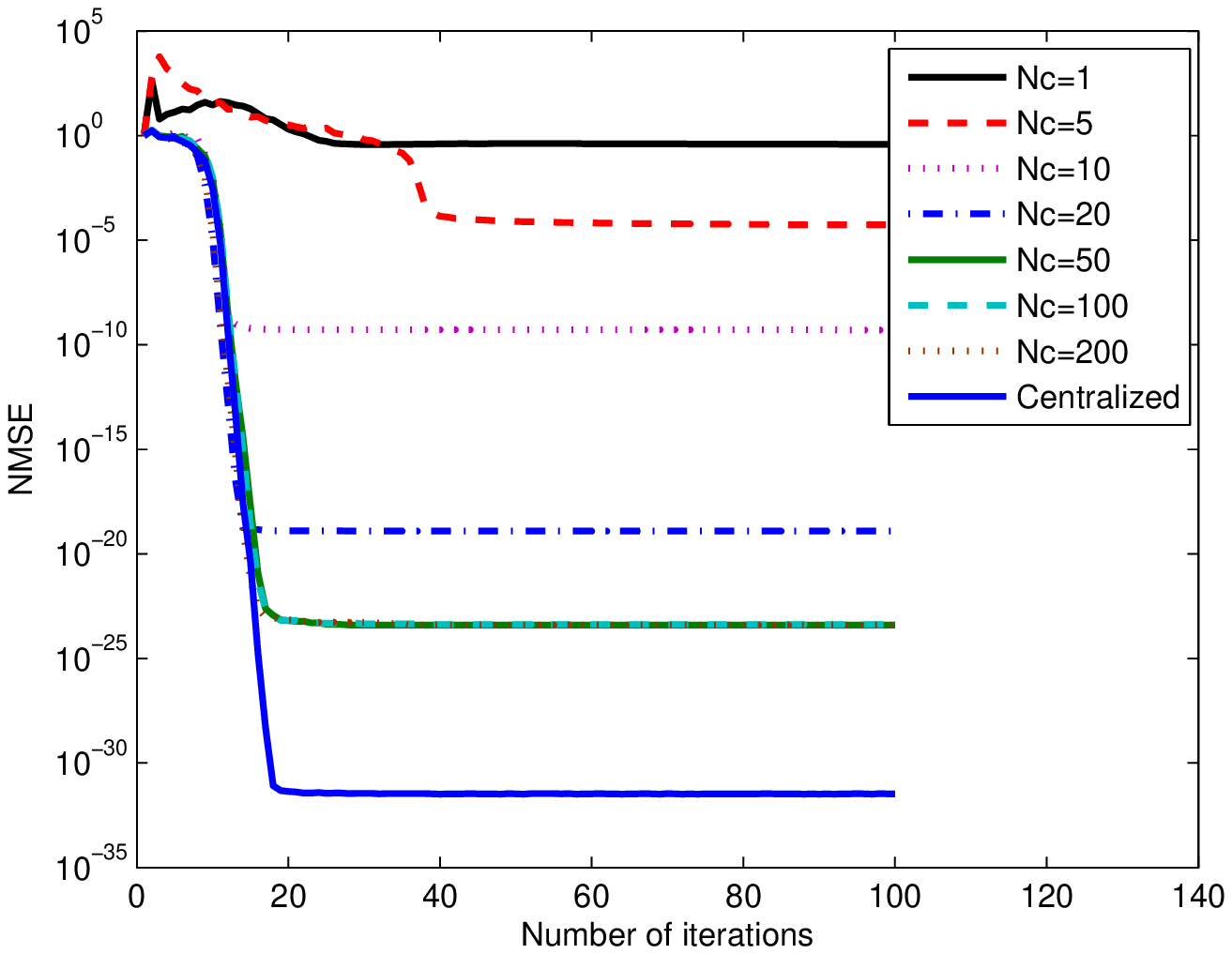}
\end{center}
\caption{Mean square for the DLM algorithm in the case of a Payley graph.}
\label{median_DLM_payley fig}
\end{figure}

Figures \ref{median_DALS_cycle_fig} and \ref{median_DLM_cycle fig} depict the NMSE for a cycle graph, while results depicted in Figures \ref{median_DALS_payley fig} and \ref{median_DLM_payley fig} concern a Payley graph. As in the centralized case, the LM based algorithms exhibit faster convergence than ALS based algorithms. However, DLM is more sensitive to imperfect consensus. Note that the performance of the DLM algorithm is improved when the number $N_c$ of consensus iterations is increased. The behavior of DALS is completely different. Moreover, the number of consensus iterations has an impact in the speed of convergence not in the final value of NMSE (see \cite{KiAl:10} for similar observations). With DALS, depending on the considered graph, even a single consensus iteration can be sufficient for achieving a good PARAFAC fitting. However, the price to pay is a slower convergence compared to DLM. We have to note that the considered cycle graph has a diameter that is twice that of the Payley graph. With a sufficient number of consensus iterations both DLM and DALS emulate their centralized counterpart.

%

\begin{figure}[!t]
\begin{center}
\includegraphics[scale=0.8]{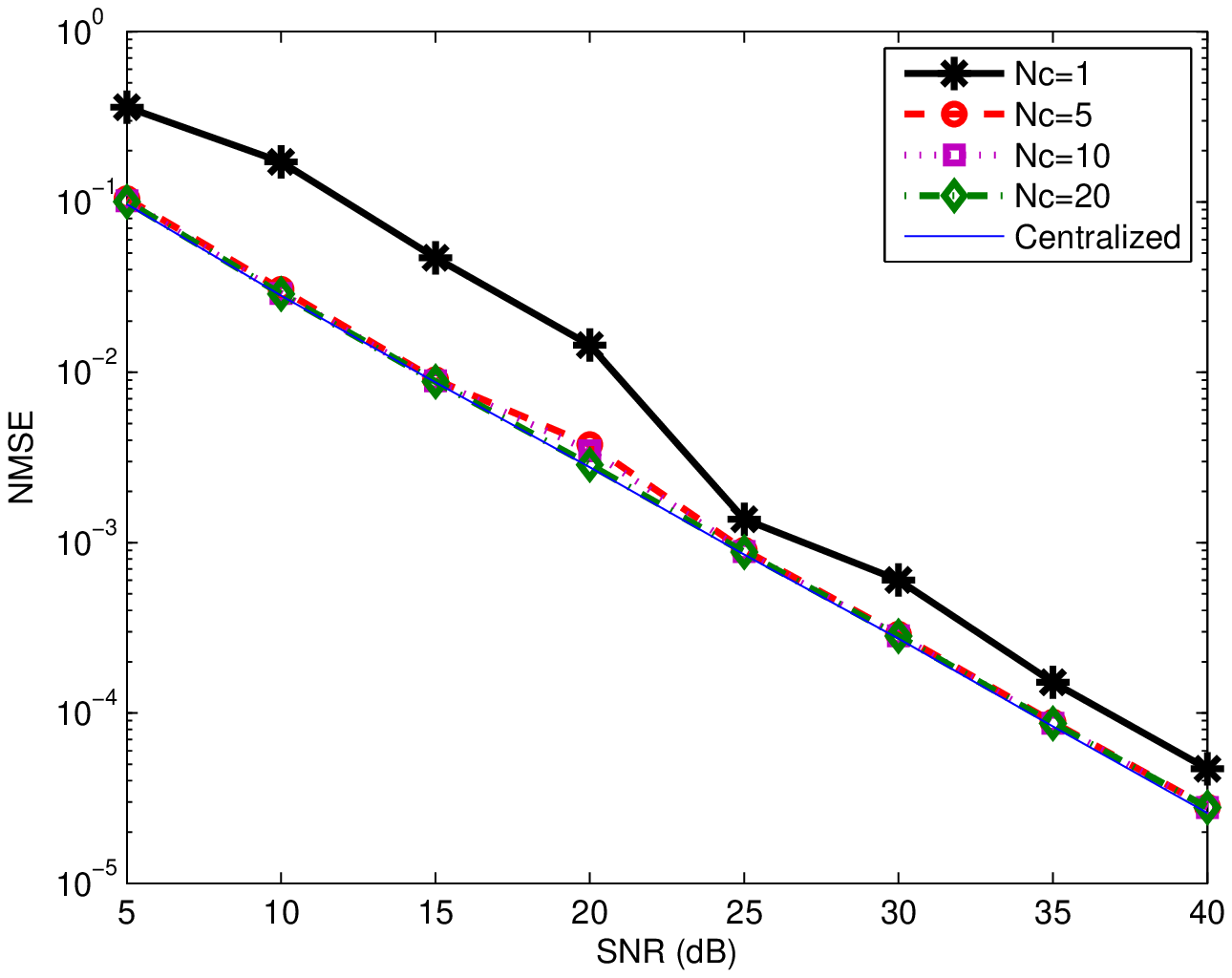}
\end{center}
\caption{Normalized Mean square error for the DALS algorithm in the case of a Payley graph (noisy case).}
\label{median_noisy_DALS_payley fig}
\end{figure}

\begin{figure}[!t]
\begin{center}
\includegraphics[scale=0.8]{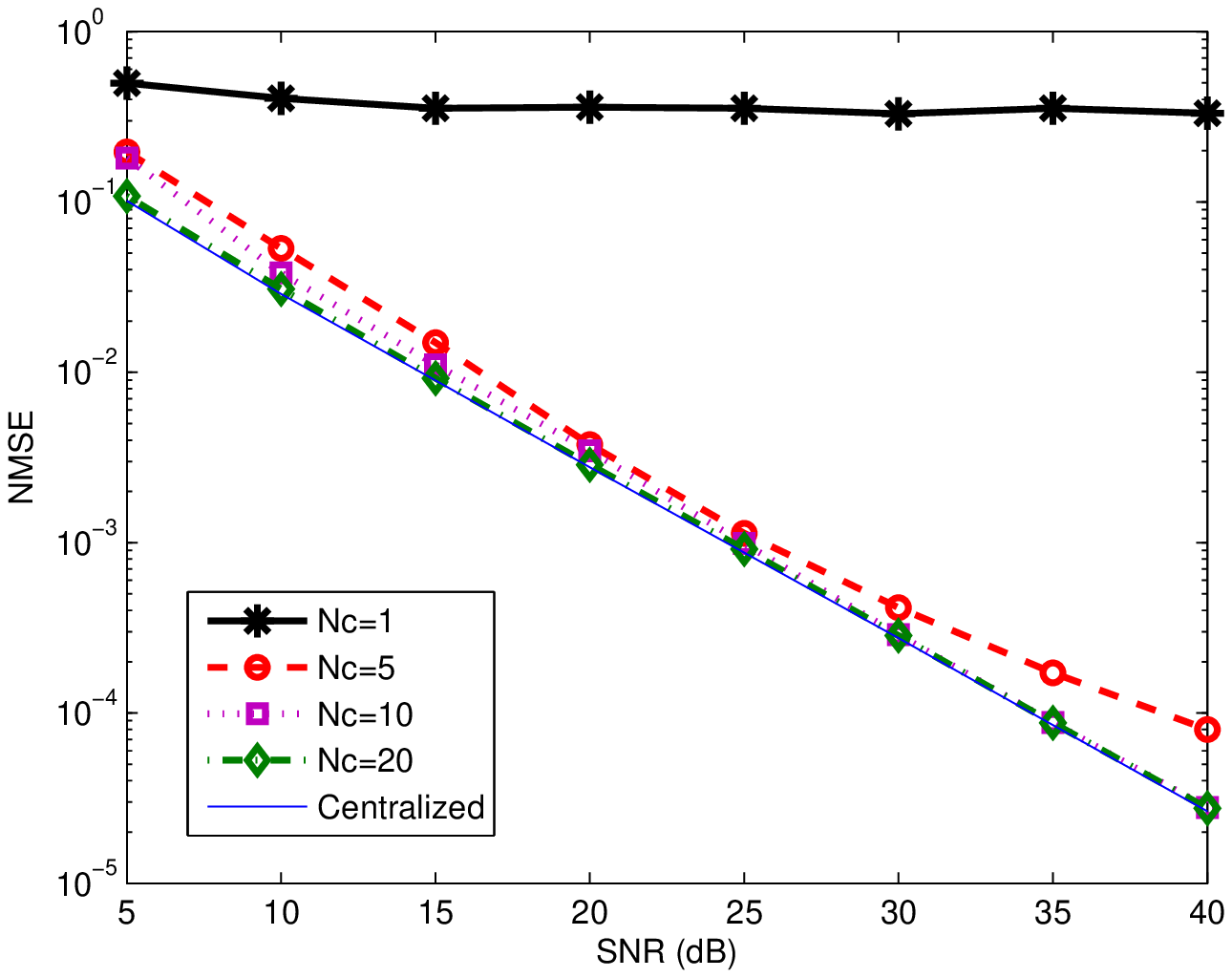}
\end{center}
\caption{Normalized Mean square error for the DLM algorithm in the case of a Payley graph (noisy case).}
\label{median_noisy_DLM_payley fig}
\end{figure}

Now, we consider noisy observations at each node. For different values of signal to noise ratio (SNR) Figures \ref{median_noisy_DALS_payley fig} and \ref{median_noisy_DLM_payley fig} depict the NMSE in steady state for DALS and DLM, respectively, when considering a Payley graph. We can note that both algorithms emulate their centralized counterparts even with a few number of iterations. As for noiseless observations, DALS emulates ALS with very few consensus iterations.

%

In another experiment, we consider a more challenging situation, where the observed data at each node are noisy and the data exchange during consensus iterations are also noisy. It is well known that standard average consensus algorithms are not robust to additive noise. In such a situation, we resort to the sequence averaging policy based on finite-time average consensus protocol proposed in \cite{Kib:12acc}. For a Payley graph, using the finite-time protocol, average consensus is reached in two steps in the noiseless case, that is the minimal time since it corresponds to the diameter of the graph. The sequence averaging protocol consists in running the finite-time protocol several time. The final result is the average of the results of the finite-time protocol. By doing so, it has been shown that exact average consensus is achieved asymptotically, i.e. when the number of sequence averaging goes towards infinity. It is worth noting that sequence averaging for noise cancelation in average consensus was first proposed in \cite{KaMo:08}. However, unlike \cite{Kib:12acc}, truncated asymptotic consensus algorithms were used.
\par
\begin{figure}[!t]
\begin{center}
\includegraphics[scale=0.8]{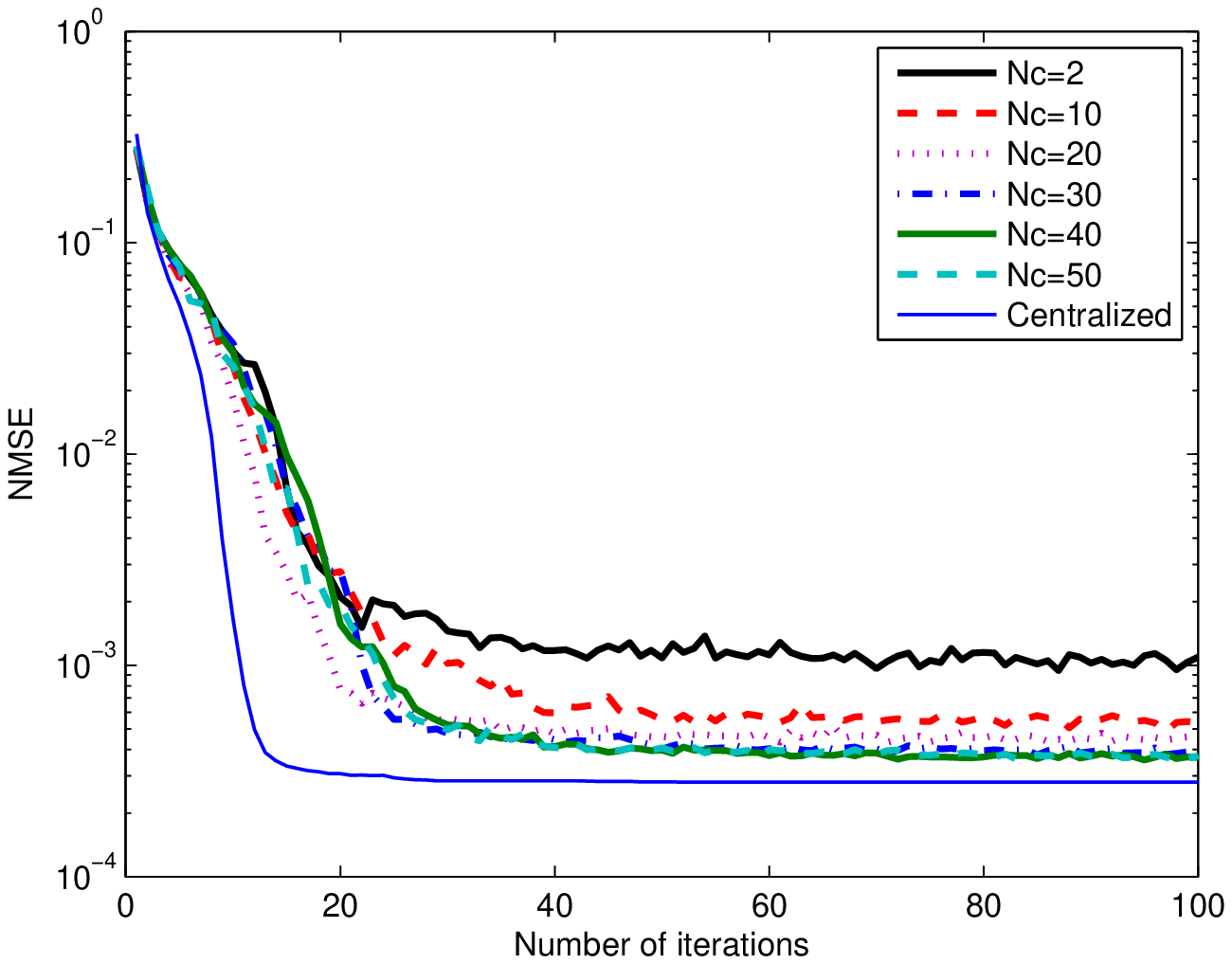}
\end{center}
\caption{Normalized Mean square error for the DALS algorithm in the case of a Payley graph with noisy data exchanges.}
\label{median_DALS_noisy_cons_payley fig}
\end{figure}

\begin{figure}[!t]
\begin{center}
\includegraphics[scale=0.8]{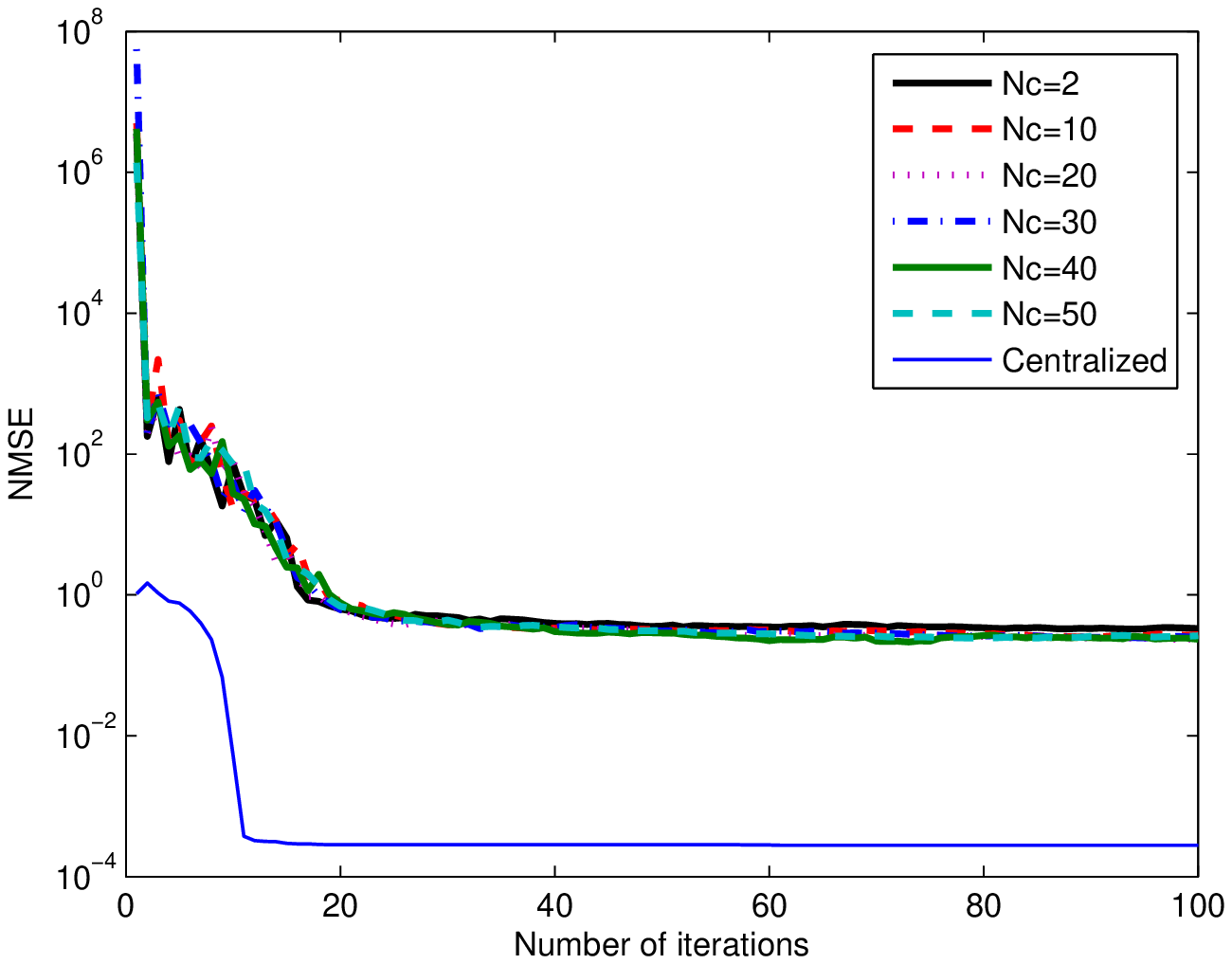}
\end{center}
\caption{Normalized Mean square error for the DLM algorithm in the case of a Payley graph with noisy data exchanges.}
\label{median_DLM_noisy_cons_payley fig}
\end{figure}

 We now compare the effect of imperfect average consensus on the performance of the proposed algorithms. Imperfections result here both from a finite number of message exchanges and noisy communications. The results depicted in Figures \ref{median_DALS_noisy_cons_payley fig} and \ref{median_DLM_noisy_cons_payley fig} were obtained with noisy observations $SNR=30 dB$ and noisy communications ($SNR=40 dB$). The number of consensus iterations $N_c=N_{ft}N_s$, with $N_s$ the number of sequence averaging and $N_{ft}$ that of finite-time average consensus protocol (here $N_{ft}=2$). We can note that DLM completely fails due to imperfect consensus while DALS exhibit a more robust behavior. By increasing the number of averaging iterations, the NMSE performance is improved and becomes closer to that of the centralized ALS algorithm.
\par From these simulations, we can conclude that the distributed algorithms proposed in this paper give similar performance than those of their centralized counterpart, under perfect average consensus. When considering imperfections in average consensus computations, DLM is more sensitive than DALS. In terms of convergence, as for the centralized case, the DLM algorithm generally has a faster convergence and yields more stable results.

\section{Conclusion}\label{sec:conc}
We have proposed a new approach to compute the PARAFAC decomposition of a third-order tensor in a distributed way across a network of collaborative nodes.
By resorting to the concept of average consensus in graphs, we have developed distributed versions of the ALS and LM algorithms that combine the estimation of local factors with an in-network computation of the global factors of the PARAFAC decomposition. Although we have restricted ourselves to the case of third-order tensors with one local and two global factor matrices for presentation simplicity, our distributed tensor-based algorithms are equally valid for $N$-th order tensors, with $N_1$ local factor matrices and $N_2$ global factor matrices, with $N=N_1+N_2$. The approach proposed in this work may be useful to a number of distributed estimation problems in signal processing, and more particularly, in the context of collaborative sensor networks in digital communications.

\bibliographystyle{IEEEtran}
\bibliography{ma_biblio}

\newpage

\end{document}